\documentclass{amsart}
\usepackage{enumitem}
\usepackage{amsmath}
\usepackage{amssymb}
\usepackage{gensymb}
\usepackage{amsthm}
\usepackage{multicol}
\usepackage[utf8]{inputenc}
\usepackage{xfrac}
\usepackage{listings}
\usepackage{xcolor}
\usepackage{pgfplots}
\usepackage{etoolbox}
\usepackage{float}
\usepackage[makeroom]{cancel}
\pgfplotsset{compat=1.17}
\usepackage{mathrsfs}
\usepackage{subcaption}

\newcommand{\spacewidth}{\the\fontdimen2\font\space}
\newcommand{\nonum}{\hspace{-\spacewidth}}

\makeatletter

\let\@newamstheorem\newtheorem

\newcommand{\@setnumber}[2]{%
  \@namedef{the#1}{#2}%
}

\newcommand{\newcustomtheorem}[2]{%
    \@newamstheorem{#1@inner}{#2}%
    \newenvironment{#1}[1][\nonum]{%
        \@setnumber{#1@inner}{##1}%
        \expandafter\begin{#1@inner}%
    }{%
        \expandafter\end{#1@inner}%
    }%
}

\renewcommand{\newtheorem}{\@ifstar{\newcustomtheorem}{\@newamstheorem}}

\makeatother
\newenvironment{alphenum}{\begin{enumerate}[label=(\alph*)]}{\end{enumerate}}

\newtheorem{thm}{Theorem}[section]
\newtheorem{lemma}[thm]{Lemma}
\newtheorem{cor}[thm]{Corollary}
\newtheorem{prop}[thm]{Proposition}
\newtheorem*{thm*}{Theorem}
\newtheorem*{lemma*}{Lemma}
\newtheorem*{prop*}{Proposition}
\newtheorem*{cor*}{Corollary}

\theoremstyle{definition}
\newtheorem{defn}[thm]{Definition}
\newtheorem*{defn*}{Definition}

\newtheorem*{ques}{Question}
\newtheorem*{prob*}{In Class Problem}

\newtheorem*{ppty*}{Property}

\newcommand{\dquote}[1]{``#1''}

\usepackage{dsfont}

\newcommand{\2}{\text{\usefont{U}{bbold}{m}{n}2}}

\newcommand{\makemathbb}[1]{%
    \expandafter\newcommand\csname #1\endcsname{\mathbb{#1}}%
}%
\newcommand{\makemathbf}[1]{%
    \expandafter\newcommand\csname b#1\endcsname{\mathbf{#1}}%
}%
\newcommand{\makemathcal}[1]{%
    \expandafter\newcommand\csname c#1\endcsname{\mathcal{#1}}%
}%
\newcommand{\makemathfrak}[1]{%
    \expandafter\newcommand\csname f#1\endcsname{\mathfrak{#1}}%
}%
\newcommand{\makemathscr}[1]{%
    \expandafter\newcommand\csname s#1\endcsname{\mathscr{#1}}%
}%

\makeatletter
\count@=0
\loop
\advance\count@ 1
\edef\upletter{\@Alph\count@}%
\edef\loletter{\@alph\count@}%
\expandafter\makemathbb\upletter
\expandafter\makemathbf\upletter
\expandafter\makemathbf\loletter
\expandafter\makemathcal\upletter
\expandafter\makemathfrak\upletter
\expandafter\makemathscr\upletter
\ifnum\count@<26
\repeat
\makeatother

\makeatletter
\def\moverlay{\mathpalette\mov@rlay}
\def\mov@rlay#1#2{\leavevmode\vtop{%
   \baselineskip\z@skip \lineskiplimit-\maxdimen
   \ialign{\hfil$\m@th#1##$\hfil\cr#2\crcr}}}
\newcommand{\charfusion}[3][\mathord]{
    #1{\ifx#1\mathop\vphantom{#2}\fi
        \mathpalette\mov@rlay{#2\cr#3}
      }
    \ifx#1\mathop\expandafter\displaylimits\fi}
\makeatother

\DeclareMathOperator{\Aut}{Aut}

\title{Shifts of Finite Type on Locally Finite Groups}
\author{Jacob Raymond}
\address{Jacob Raymond\\
Department of Mathematics and Statistics\\
University of North Carolina at Charlotte \\
9201 University City Blvd.\\
Charlotte, NC 28223}
\email{jraymon9@uncc.edu}

\begin{document}

\begin{abstract}
    In this work, we prove that every SFT, sofic shift, and strongly irreducible shift on locally finite groups has strong dynamical properties. These properties include that every sofic shift is an SFT, every SFT is strongly irreducible, every strongly irreducible shift is an SFT, every SFT is entropy minimal, and every SFT has a unique measure of maximal entropy, among others. In addition, we show that if every SFT on a group is strongly irreducible, or if every sofic shift is an SFT, then the group must be locally finite, and this extends to all of the properties we explore. These results are collected in two main theorems which characterizes the local finiteness of groups by purely dynamical properties. In pursuit of these results, we present a formal construction of \textit{free extension} shifts on a group $G$, which takes a shift on a subgroup $H$ of $G$, and naturally extends it to a shift on all of $G$. 
\end{abstract}

\maketitle

\section{Introduction} \label{INT}

For a finite set of symbols $\cA$ and a group $G$, the field of symbolic dynamics studies the action of $G$ by translations on the set $\cA^G$, called \textit{full $G$-shift with alphabet $\cA$}, and the subsystems within. Equipped with the product topology (with the discrete topology on $\cA$), a closed, translation invariant subset of $\cA^G$ is called a \textit{$G$-shift}, and understanding what properties such subsystems can exhibit is central to symbolic dynamics. In its conception, the primary group of interest was $\Z$, the group of integers under addition. Even in this case, complex behavior arises, though much is known in general about shifts on $\Z$ \cite{lind_marcus}. A natural extension of this case is the group $\Z^d$ for some natural number $d$, the study of which has been called multi-dimensional symbolic dynamics. More recently, interest in shifts on $\Z^d$ has grown, though this case already adds much complexity \cite{robinson, mcgoff_random, hochman}, and less is known about $\Z^d$-shifts in general. Interest in the general group case is even more recent, and as may be expected, is even less tractable than the case of $\Z^d$, though a recent result about tilings of amenable groups \cite{dhz} has made a few results about shifts on amenable groups possible \cite{frisch_tamuz, bland_ss_ent, bland_embedding}.

The class of \textit{$G$-shifts of finite type}, or $G$-SFTs, are of particular interest, as they are characterized by a finite amount of information. More precisely, a $G$-SFT $X$ is a $G$-shift for which there is a finite collection of \textit{patterns} (an element of $\cA^F$ for a finite $F \subset G$) so that $X$ is the collection of all \textit{configurations} in $\cA^G$ for which these patterns never appear. The finite nature of $G$-SFTs makes them amenable to analysis using finitary and combinatorial methods, and in general $G$-SFTs are well behaved in comparison to general shifts. Furthermore, every shift on a group can be represented as an intersection of SFTs, so in this sense, SFTs are plentiful and are good approximations for shifts in general. Formal definitions of $G$-shifts and $G$-SFTs can be found in Section \ref{DEF:shifts}.

Understanding what properties are possible for SFTs on groups is at the core of symbolic dynamics. One such property is the \textit{entropy} (Definition \ref{DEF:top_ent}) of an SFT on a countable amenable group $G$, or in particular, the set of entropies which are attainable by SFTs on $G$, which is denoted $\cE(G)$. $\cE(\Z)$ was classified by Lind \cite{lind_Z_entropies}, and more recently, $\cE(\Z^d)$ for $d \ge 2$ was classified by Hochman and Meyerovitch \cite{hochman_meyerovitch}. Recent results by Barbieri \cite{barbieri} classifies $\cE(G)$ as $\cE(\Z^d)$ for a certain class of amenable groups. Currently to the knowledge of the author, there are no known finitely generated groups $G$ for which $\cE(G)$ does not coincide with either $\cE(\Z)$ or $\cE(\Z^2)$, and further classifying $\cE(G)$ for other groups and classes of groups is an open goal in symbolic dynamics. Another property is \textit{strong irreducibility} (Definition \ref{DEF:si}), which loosely gives that a $G$-shift is large, and contains a large variety of configurations. In general, a $G$-SFT need not be strongly irreducible, and a strongly irreducible $G$-shift need not be a $G$-SFT. The additional structure which strong irreducibility imposes on a shift has been useful in proving results about shifts \cite{ceccherini_silberstein, pavlov, bland_embedding}. We also explore several other properties of shifts, which are outlined in Section \ref{DEF}, and discussed informally after the statement of our two main theorems below.

Our motivation for studying locally finite groups comes from the following example. Let $G = \bigoplus_{n \in \N} \Z / 2 \Z$, the countable direct sum of the 2 element group. Elements of $G$ are infinite sequences of 0s and 1s which only contain finitely many 1s, and the group operation is component-wise addition modulo 2. Using elementary methods for computing the entropy on shifts, it is possible to show directly that
\[
    \cE(G) = \left\{\frac{\log(n)}{2^{m-1}} : n,m \in \N\right\} \subsetneq \cE(\Z),
\]
providing an example of an infinitely generated group for which $\cE(G)$ does not coincide with $\cE(\Z)$ or $\cE(\Z^2)$. In general, classifying the entropies which are attainable by SFTs for a group $G$ is quite difficult, however the process is made tractable for this group by the fact that
\[
    H_n = \left(\bigoplus_{k=1}^n \Z / 2\Z\right) \oplus \left(\bigoplus_{k=n+1}^\infty \{0\}\right)
\]
is a sequence of finite subgroups of $G$ such that $H_n \le H_{n+1}$ and $G = \bigcup_{n \in \N} H_n$, which makes $\{H_n\}$ a F{\o}lner sequence for $G$. As it turns out, a countable group with such a sequence $\{H_n\}$ of finite subgroups is necessarily \textit{locally finite}. A group is locally finite if every finitely generated subgroup is finite. In fact, any countable locally finite group must have such a sequence of subgroups, and so this property coincides exactly with the property of being locally finite when the group is countable. Locally finite groups naturally extend finite groups in a way that allows for finitary methods to be used when analyzing the groups, despite being possibly infinite. As a result, one may suspect SFTs on locally finite groups are highly structured and have many nice dynamical properties.

The main results of this paper confirm that SFTs on locally finite groups have very strong dynamical properties. Furthermore, we show that locally finite groups are the only groups for which all SFTs exhibit these properties. These results are grouped in two, one in the case where $G$ is an arbitrary group, and the second where $G$ is a countable amenable group. The first is given below, and followed by a brief explanation of each statement in the result, though formal definitions for every term below can be found in Section \ref{DEF}. 

\begin{thm*}[I] \label{thmI}
    Let $G$ be a group. Then the following are equivalent.
    \begin{alphenum}
        \item \label{I:lf} $G$ is locally finite.
        \item \label{I:fe} Every $G$-SFT is the free extension of some SFT on a finite subgroup of $G$.
        \item \label{I:ir} Every $G$-SFT is strongly irreducible.
        \item \label{I:ir2} Every strongly irreducible $G$-shift is a $G$-SFT.
        \item \label{I:sof} Every sofic $G$-shift is a $G$-SFT.
        \item \label{I:aut} For every $G$-SFT $X$, $\Aut(X)$ is locally finite.
    \end{alphenum}
\end{thm*}

Statement \ref{thmI}\ref{I:fe} is not a typical dynamical property, but involves a specific type of shift defined in Section \ref{FE} called a \textit{free extension} shift. Free extensions shifts are by no means a new concept and have been used in the past \cite{hochman_meyerovitch, barbieri}, however we present a formal construction which proves to be quite useful in deriving properties of free extensions. The equivalence between statement \ref{thmI}\ref{I:fe} and \ref{thmI}\ref{I:lf} is at the core of nearly every argument involved in proving this theorem and the next. Free extension shifts are defined for general groups in Section \ref{FE}, and may be useful in studying shifts on groups in general, beyond the study of shifts on locally finite group. Statement \ref{thmI}\ref{I:sof} involves sofic shifts, which are the image of SFTs under continuous, shift invariant factor maps. Along with SFTs, sofic shifts are a noteworthy class of shifts which are defined by a finite amount of information. Every SFT is necessarily sofic, however the converse does not hold in general, and Theorem \ref{thmI} gives that the converse holds only in the case that the group is locally finite. The definition of factor maps can be found in Section \ref{DEF:factors}, and sofic shifts in Section \ref{DEF:sofic}. Statement \ref{thmI}\ref{I:ir} gives that every SFT on a locally fiinte group is strongly irreducibile. A formal definition is given by Definition \ref{DEF:si}, but informally, strong irreducibility is a property which guarantees that for any two elements in the shift, there exists an element of the shift which is equal to one of the elements on a finite subset, and equal to the other on any sufficiently separated finite subset. In this sense, strongly irreducible shifts are rich with configurations. Statement \ref{thmI}\ref{I:ir2} is the converse of the previous statement, and is independently equivalent to the group being locally finite. These two statements in combination give that the set of $G$-SFTs and the set of strongly irreducible $G$-shifts coincides exactly when $G$ is locally finite, but that neither is contained in the other when $G$ is not locally finite. Statement \ref{thmI}\ref{I:aut} involves $\Aut(X)$, the automorphism group of an SFT $X$. This group consists of homeomorphisms from $X$ to itself which preserve the action of $G$, and is formally defined at the end of Section \ref{DEF:factors}.

For the second result, we restrict to the case that $G$ is a countable amenable group, which permits the development of topological entropy, and each of the statements in the result involves this entropy. A brief discussion of each statement follows the statement of the result, and the formal definitions of every term can be found in Section \ref{DEF}. In statement \ref{thmII}\ref{II:sft_ent}, we use the non-standard notation $H \ll G$ to denote that $H$ is a \textit{finite} subgroup of $G$.

\begin{thm*}[II] \label{thmII}
    Let $G$ be a countable amenable group. Then the following are equivalent.
    \begin{alphenum}
        \item \label{II:lf} $G$ is locally finite 
        \item \label{II:0_ent} If $X$ is a nonempty $G$-SFT with $h(X) = 0$, then $X = \{x\}$, where $x$ is a fixed point.
        \item \label{II:ent_min} Every $G$-SFT is entropy minimal.
        \item \label{II:sft_ent} $G$ is locally typical and
        \[
            \cE(G) = \left\{\frac{\log(n)}{|H|} : H \ll G, n \in \N \right\} \subset \Q^+_{\log} = \left\{\frac{\log(n)}{m} : n,m \in \N\right\}.
        \]
        \item \label{II:mme} Every $G$-SFT has a unique measure of maximal entropy.
    \end{alphenum}
\end{thm*}

We remark that while we restrict the results to countable amenable groups, entropy can also be extended to the more general class of countable sofic groups \cite{bowen_sofic}, however we will not need this more general definition, since any countable locally finite group is necessarily amenable. The definition of entropy can be found in Section \ref{DEF:ent}. Statement \ref{thmII}\ref{II:0_ent} is about what sorts of zero topological entropy SFTs can exist, and in the case of locally finite groups, there is a single zero-entropy SFT (up to conjugacy). This result indirectly answers a question of Barbieri in the affirmative:

\begin{ques}[3.19 \cite{barbieri}]
    Does there exist an amenable group $G$ and a $G$-SFT which does not contain a zero-entropy $G$-SFT?
\end{ques}

Since the only 0 entropy SFTs on locally finite groups are single fixed points, it suffices to construct an SFT which contains no fixed points, which is trivial to do using free extensions. There is further discussion about this construction in Section \ref{FR}.

Statement \ref{thmII}\ref{II:ent_min} involves entropy minimality, which is the property that a shift has no proper subshift with the same entropy as the entire shift. A formal definition is given by Definition \ref{DEF:ent_min}. Statement \ref{thmII}\ref{II:sft_ent} consists of two parts. The first, is that $G$ is locally typical, which is non-standard terminology for groups, and means that every finitely generated subgroup of $G$ is either finite, or contains an element of infinite order. The need for this requirement in this statement is discussed further in Sections \ref{LF-sft_ent} and \ref{FR}. The second part of the statement classifies the set of entropies attainable on any locally finite group, and gives the following corollary which may be of independent interest to the remainder of the Theorem.

\begin{cor*}
    Let $G$ be a countable locally finite group. Then
    \[
        \cE(G) = \left\{\frac{\log(n)}{|H|} : H \ll G, n \in \N\right\}
    \]
\end{cor*}

Finally, the last statement in the theorem, Statement \ref{thmII}\ref{II:mme}, involves measure theoretic entropy and measures of maximal entropy, which are invariant measures on the SFT that have a measure theoretic entropy equal to the topological entropy of the system. Formal definitions for these can be found starting at Definition \ref{DEF:mes_ent}.

\subsection{Overview}

In Section \ref{DEF}, we present the relevant background and notation used in the remainder of the paper. In Section \ref{FE}, we construct free extension shifts generally for groups, and then prove some properties of these shifts. In Section \ref{LF} we prove Theorems \ref{thmI} and \ref{thmII}, which is broken down into several individual lemmas. Finally, in Section \ref{FR}, we discuss some general consequences of Theorems \ref{thmI} and \ref{thmII} and properties of free extensions, and indicate possible directions for future work.

\section{Definitions and Notation} \label{DEF}

We begin with defining all necessary background terms and notation. The section is broken up into subsections based on what is being defined.

\subsection{Sets and Groups}

For any set $A$, let $B \Subset A$ denote that $B$ is a \textit{finite} subset of $A$. The set difference of two sets $A$ and $B$ is denoted by $A \setminus B$. The disjoint union of two sets $A$ and $B$ is denoted by $A \sqcup B$. The symmetric difference of two sets $A$ and $B$ is denoted by $A \triangle B$

Given two sets $A$ and $B$, the set $A^B$ refers to the collection of all functions $f:B \to A$. If $A$ is endowed with a topology, then $A^B$ is endowed with the product topology.

For a group $G$, we denote that a subset $H \subset G$ is a subgroup of $G$ by $H \le G$, and to additionally specify that $H$ is a finite subgroup of $G$, we use the notation $H \ll G$. For $F \subset G$, the subgroup of $G$ generated by $F$ is denoted as $\langle F \rangle$. A group is \textit{periodic} (also known as a \textit{torsion group}) if all of its elements have finite order. These notations are standard for groups, however the following is non-standard; A group is called \textit{atypical} if it is periodic and infinite. A group is called \textit{typical}, if it is not atypical, which means it is finite or not periodic. If $P$ is a property which a group can posses, then a group $G$ is said to be \textit{locally $P$} if $\forall F \Subset G$, the subgroup $\langle F \rangle$ has property $P$. A group $G$ is then \textit{locally finite} if $\forall F \Subset G$ we have  $\langle F \rangle \ll G$, that is, every finitely generated  subgroup of $G$ is finite. $G$ is \textit{locally typical} if $\forall F \Subset G$, the subgroup $\langle F \rangle$ is typical. In addition, we use the terminology \textit{non-locally finite} to mean that a group is not locally finite (and similarly for non-locally typical).

Given a group $G$ and subgroup $H \le G$, we denote the set of \textit{right} cosets of $H$ in $G$ by $H \backslash G$. This notation is similar to the one used for set difference (though the spacing is different), however which is being referred to is generally clear from context.

A countable group $G$ is \textit{amenable} if there exists a sequence $\{F_n\}_{n = 1}^\infty$ such that $F_n \Subset G$, $\{F_n\}$ exhausts $G$ so that $G = \bigcup_n F_n$, and $\forall g \in G$,
\[
    \lim_{n \to \infty} \frac{|gF_n \triangle F_n|}{|F_n|} = 0.
\]
Such a sequence is called a \textit{left F{\o}lner sequence}, and similarly, \textit{right F{\o}lner sequences} exist for amenable groups, which satisfy
\[
    \lim_{n \to \infty} \frac{|F_ng \triangle F_n|}{|F_n|} = 0.
\]

\subsection{$G$-spaces}

Though the main results of this paper do not pertain to general $G$-spaces, they are used in constructions, and are therefore necessary to define. $G$-spaces are defined more generally for topological groups, but in the particular case of this paper, the topology on the group $G$ is always taken to be the discrete topology, so the definition reduces to the following.

\begin{defn}
    Let $G$ be a group (viewed as a discrete group), and $X$ be a topological space. Take any collection $T = \{T^g : g \in G\}$ of functions $T^g : X \to X$ which are homeomorphisms, and satisfy $T^g \circ T^h = T^{gh}$ for all $g,h \in G$. Additionally, it should be that $T^e$ is the identity function on $X$ for the identity element $e \in G$. Then the pair $(X, T)$ is a $G$-space.
\end{defn}

In our case, the topological space $X$ is always a Hausdorff space, usually compact, and occasionally metrizable.

\subsubsection{Factors and Conjugacies} \label{DEF:factors}

The main use for general $G$-spaces in this paper comes in the form of conjugacies between two $G$-spaces, a natural notion of equivalence. First, we define factors and factor maps.

\begin{defn}
    Let $G$ be a group, and let $(X, T)$, $(Y, S)$ be two $G$-spaces. Then a map $\phi:X \to Y$ is a \textit{factor map} if
    \begin{itemize}
        \item $\phi$ is continuous,
        \item $\phi$ is surjective, and
        \item for every $g \in G$, $S^g \circ \phi = \phi \circ T^g$.
    \end{itemize}
    If such a factor map exists from $X$ to $Y$, then $Y$ is called a \textit{factor} of $X$.
\end{defn}

In the case that a factor map $\phi$ is a homeomorphism, then $\phi$ is called a \textit{conjugacy}, and $X$ and $Y$ are said to be \textit{conjugate}. The collection of conjugacies from a $G$-space $X$ to itself forms a group under composition, denoted $\Aut(X)$.

\subsection{$G$-shifts} \label{DEF:shifts}

For the remainder of the section, $\cA$ is a finite \textit{alphabet} (set), endowed with the discrete topology.

\begin{defn}
    Let $G$ be a group. $\cA^G$ is endowed with the product topology, which makes it a compact Hausdorff space. When $G$ is countable $\cA^G$ is metrizable, and we will take this fact to be evident when $G$ is countable. For any $g \in G$, define $\sigma^g:\cA^G \to \cA^G$ by
    \[
        (\sigma^g x)(h) = x(hg)
    \]
    for any $h \in G$. Each $\sigma^g$ is a homeomorphism from $\cA^G$ to $\cA^G$, and $\sigma^g \circ \sigma^h = \sigma^{gh}$ for all $g, h \in G$. Also, $\sigma^ex = x$ for each $x \in X$, and therefore the collection $\sigma = \{\sigma^g : g \in G\}$ is a continuous action of $G$ on $\cA^G$. The pair $(\cA^G, \sigma)$ is called the \textit{full $G$-shift with alphabet $\cA$}, or simply the \textit{full $G$-shift} when the alphabet $\cA$ is clear, which is typically the case. This also makes $(\cA^G, \sigma)$ a $G$-space. The elements of $\cA^G$ are referred to as \textit{configurations}.
\end{defn}

Though the full $G$-shift is interesting in its own right, we are primarily interested in subsystems of the full $G$-shift, which are called \textit{$G$-shifts}.

\begin{defn}
    Let $G$ be a group. A subset $X \subset \cA^G$ is said to be \textit{$G$-invariant}, or merely \textit{shift invariant} when the group $G$ is clear from context, if for every $x \in X$, and $g \in G$, $\sigma^g x \in X$. A closed, $G$-invariant subset $X \subset \cA^G$, along with the action of $G$ on $\cA^G$ restricted to $X$, is called a \textit{$G$-shift of $\cA^G$}, or just \textit{$G$-shift} when the full shift is clear from context. $G$-shifts are also $G$-spaces.
\end{defn}

\subsubsection{Patterns}

Although an element of $\cA^G$ is known as a configuration, the term \textit{pattern} is used when considering elements of $\cA^F$ for some $F \subset G$. In addition, we define a few operations on patterns that are quite useful when working with shifts.

\begin{defn}
    For any group $G$ and $F \subset G$, an element $w \in \cA^F$ is called an \textit{$\cA$-pattern on $F$}, or just a \textit{pattern} if $\cA$ is clear. The \textit{shape} of a pattern $w \in \cA^F$ is the set $F$ itself.
    
    For $E, F \subset G$ and patterns $w \in \cA^E$ and $v \in \cA^F$, we say $w$ and $v$ are disjoint if $E$ and $F$ are disjoint. Similarly, $w \in \cA^F$ is said to be finite if $F$ is finite, and infinite if $F$ is infinite.
    
    For any $E \subset F \subset G$ (including $E = F = G$), the restriction of a pattern $w \in \cA^F$ to $E$, which is denoted $w|_E$ and contained in $\cA^E$, is defined as $w|_E(g) = w(g)$ for every $g \in E$. Conversely, for some $w \in \cA^E$, the set of \textit{$F$-extensions of $w$} is defined as
    \[
        [w]_F = \bigl\{v \in \cA^F : v|_E = w\bigr\}.
    \]
    In the case that $F \le G$, then $[w]_F$ is known as a \textit{cylinder set}. In the case that $F = G$, then $[w]$ is used instead of $[w]_G$, unless clarity is necessary.
\end{defn}

Patterns are very useful in describing the structure of $G$-shifts. For any $G$-shift $X$ (including the full $G$-shift), the set
\[
    \fB = \bigl\{[w]_G \cap X : F \Subset G, w \in \cA^F\bigr\}
\]
is a basis for the subspace topology on $X$ as a subspace of the full $G$-shift. Note that $[w]_G \cap X$ may be empty or nonempty, and we define the following sets in order to distinguish when this is or is not the case.

\begin{defn}
    For any $G$-shift $X$, and any $F \subset G$, let $\cL_F(X)$ denote the \textit{$F$-language} of $X$, which is defined as
    \[
        \cL_F(X) = \{x|_F : x \in X\} \subset \cA^F.
    \]
    We then let $\cL(X)$ be the \textit{language} of $X$, which is defined as
    \[
        \cL(X) = \bigsqcup_{F \Subset G} \cL_F(X).
    \]
    By this definition, note that $w \in \cL(X)$ if and only if $[w]_G \cap X \ne \varnothing$. In addition, let $\cL^\infty(X)$ denote the set
    \[
        \cL^\infty(X) = \bigsqcup_{F \subset G} \cL_F(X).
    \]
    The main difference between this and $\cL(X)$ is that $\cL^\infty(X)$ also contains infinite patterns.
    
    We also let $\cF_F(X) = \cA^F \setminus \cL_F(X)$, and 
    \[
        \cF(X) = \bigsqcup_{F \Subset G} \cF_F(X).
    \]
    These sets are known as the \textit{forbidden $F$-patterns} of $X$ and the \textit{forbidden patterns} of $X$, respectively.
\end{defn}

In constructions which appear in Section \ref{FE}, we utilize an extension of the shift action $\sigma$ to $\cL^{\infty}(X)$, as well as a joining operation which allows taking two disjoint patterns and combining them into one pattern. These are defined next.

\begin{defn}
    Let $G$ be a group, and $X$ be a $G$-shift. Let $g \in G$. Then for any $F \subset G$, define $\sigma_F^g:\cL_F(X) \to \cL_{Fg^{-1}}(X)$ by
    \[
        (\sigma_F^gw)(h) = w(hg), \quad \forall h \in Fg^{-1}.
    \]
    Note that in the case $F = G$, this covers the typical shift maps. We then define $\sigma^g:\cL^\infty(X) \to \cL^\infty(X)$ for any $F \subset G$ and pattern $w \in \cA^F$ as
    \[
        \sigma^gw = \sigma_F^gw.
    \]
    In fact, $(\cL^\infty(X), \sigma)$ is a $G$-space when $\cL^\infty(X)$ is endowed with the disjoint union topology (where $\cL_F(X)$ is endowed with the subspace topology from the product topology on $\cA^F$).

    Restricting patterns to subshapes and shifting behave well in relation to each other. Let $E \subset F \subset G$, and let $g \in G$. Then for any $w \in \cL_{Eg}(X)$, the pattern $\sigma^gw$ has shape $Egg^{-1} = E$, and for any $h \in E$,
    \[
        \sigma^g(w|_{Eg})(h) = (w|_{Eg})(hg) = w(hg) = (\sigma^gw)(h).
    \]
    Since this holds for any $h \in E$, it follows that
    \[
        \sigma^g(w|_{Eg}) = (\sigma^gw)|_E.
    \]
    This rule is used in many proofs without reference.
    
    Similar interplay exists between the shifts and extension sets. Let $E \subset F \subset G$, and $g \in G$. Then $Eg \subset Fg$, and for any $w \in \cL_{Eg}(X)$,
    \[
        \sigma^g[w]_{Fg} = [\sigma^gw]_F.
    \]
    This is also used in many proofs without reference.
\end{defn}

Along with this natural notion of shifting patterns, there is a natural way to define joining two disjoint patterns.

\begin{defn}
    Let $G$ be a group, and $X$ be a $G$-shift. For any disjoint $u, v \in \cL^\infty(X)$, with shapes $F_u$ and $F_v$ respectively (so that $F_u \cap F_v = \varnothing$), we define the \textit{join} of $u$ and $v$, denoted by $u \vee v$, as follows. Let $w = u \vee v$ be defined as
    \[
        w(g) = \begin{cases}
            u(g), & g \in F_u \\
            v(g), & g \in F_v
        \end{cases},
    \]
    which is a pattern with shape $F_u \sqcup F_v$. Since $F_u$ and $F_v$ must be disjoint to take a join, it is clear that $\vee$ is commutative.
    
    Additionally, the shift action distributes over $\vee$. For any disjoint $u, v \in \cL^\infty(X)$ and $g \in G$, it is always the case that $\sigma^g(u \vee v) = (\sigma^gu) \vee (\sigma^gv)$.
    
    Furthermore, for any infinite collection of mutually disjoint patterns, all of these patterns can be joined together into one (possibly infinite) pattern, and by this commutativity, the order of the infinite join is irrelevant. Also, the shifts commute with infinite joins for similar reasons. Infinite joins and the commutativity of the shifts with infinite joins are an integral part of Section \ref{FE}.
\end{defn}

\subsubsection{Properties of $G$-shifts}

Each $G$-shift $X$ defines a set of forbidden patterns, however it is also possible to define a $G$-shift from a set of forbidden patterns.

\begin{defn}
    Let $G$ be a group, $\cA$ be a finite alphabet, and let $\bF \subset \cL(\cA^G)$ be a set of forbidden patterns. Define
    \[
        \cX^G[\bF] = \bigl\{x \in \cA^G : \forall g \in G, \forall F \subset G, \: (\sigma^gx)|_F \not\in \bF\bigr\}
    \]
    It is an elementary exercise to show that $\cX^G[\bF]$ is a $G$-shift (though possibly empty), so $\cX^G[\bF]$ is called the \textit{$G$-shift defined by $\bF$}. $\cX[\bF]$ is used whenever $G$ is clear from the context.
    
    Another elementary result is that $X = \cX[\cF(X)]$ for any $G$-shift $X$, and therefore every $G$-shift is generated by some set of finite forbidden patterns.
\end{defn}

While $\cF(X)$ is always a set of forbidden patterns which defines the $G$-shift $X$, there may be much smaller sets of forbidden patterns which exist. In some cases, there may be a finite set of forbidden patterns which defines a $G$-shift $X$, in which case the $G$-shift is called a \textit{shift of finite type}.

\begin{defn}
    Let $G$ be a group, and $X$ a $G$-shift. Then $X$ is called a \textit{$G$-shift of finite type}, or typically a \textit{$G$-SFT}, if there exists a finite $\bF \Subset \cL(X)$ such that $X = \cX[\bF]$.
    
    For any $G$-SFT $X$ there always exists some $F \Subset G$ such that $X = \cX[\cF_F(X)]$. Such a shape $F$ is called a \textit{forbidden shape} for $X$. Additionally, given some forbidden shape $F$, any $H \Subset G$ with $F \subset H$ is also a forbidden shape, meaning $X = \cX[\cF_F(X)] = \cX[\cF_H(X)]$. This property is used in many results without reference.
\end{defn}

The finitary nature of $G$-SFTs makes them amenable to analysis using more combinatorial methods, and they are generally well behaved in many regards. Another strong property a $G$-shift can possess is \textit{strong irreducibility}, which is a strong mixing type property that is of general interest in the literature. 

\begin{defn} \label{DEF:si}
    Let $G$ be a group, and $X$ be a $G$-shift. Then $X$ is \textit{strongly irreducible} if there exists a finite $K \Subset G$ with the following property. For any $u, v \in \cL(X)$ with shapes $F_u$ and $F_v$, if $F_u \cap KF_v = \varnothing$, then there exists $x \in X$ such that $x|_{F_u} = u$ and $x|_{F_v} = v$.
    
    This definition differs from typical definitions of strong irreducibility of shifts on finitely generated groups \cite{frisch_tamuz}. In the case that $G$ is finitely generated, this definition is equivalent to more typical definitions, and is merely an extension of the more typical definition to (possibly) infinitely generated groups.
\end{defn}

\subsubsection{Factors and Sofic shifts} \label{DEF:sofic}

While for general $G$-spaces, not much can be said about the structure of factor maps, for $G$-shifts (and a broader class of $G$-shift-like $G$-spaces), factor maps have a very specific structure. We begin by defining a specific kind of factor map which can be constructed between two $G$-shifts.

\begin{defn}
    Let $\cA$ and $\cB$ be two finite alphabets, let $G$ be a group, and let $X$ be a $G$-shift of $\cA^G$. For some $F \Subset G$, let $\beta:\cL_F(X) \to \cB$ be any function, called a \textit{block map}. Then $\beta$ induces a map $\phi^G_\beta:X \to \cB^G$ called a \textit{block code} by
    \[
        \bigl(\phi^G_\beta(x)\bigr)(g) = \beta\bigl((\sigma^g x)|_F\bigr),
    \]
    and $Y = \phi^G_\beta(X)$ is a $G$-shift of $\cB^G$. Rather than $\cB^G$ however, we consider the co-domain of $\phi^G_\beta$ to be $Y$, which makes $\phi^G_\beta$ surjective and therefore a factor map from $X$ to $Y$.
\end{defn}

Block codes are generally easy to work with, due to the finitary nature of the block map that generates them. Surprisingly, any factor map between $G$-shifts (on possibly different alphabets) is a block code generated by some block map, and this fact is given by the following theorem.

\begin{thm}[Curtis-Lyndon-Hedlund] \label{clh}
    Let $G$ be a group, $\cA$ and $\cB$ be finite alphabets, $X$ be a $G$-shift of $\cA^G$ and $Y$ a $G$-shift of $\cB^G$, and let $\phi:X \to Y$. Then $\phi$ is a factor map if and only if there exists $F \Subset G$ and block map $\beta:\cL_F(X) \to \cB$ such that $\phi = \phi^G_\beta$.
\end{thm}

A proof of the theorem at this level of generality can be found in \cite[Corallary 6]{chl}. Informally, the theorem gives that factor maps for $G$-shifts are defined by a finite amount of information. A broader class of $G$-shifts which are defined by a finite amount of information, which contains all SFTs but generally includes more shifts, is the class of sofic $G$-shifts.

\begin{defn}
    A $G$-shift $Y$ is called a \textit{sofic $G$-shift} if there exists a $G$-SFT $X$ such that $Y$ is a factor of $X$.
\end{defn}

Weiss noted when first introducing sofic $\Z$-shifts that \dquote{the finite type subshifts are flawed by not being closed under the simplest operation, namely that of taking [factors]} \cite{weiss_sofic}. The collection of all sofic shifts is clearly closed under taking factors, and this is one of the many reasons the class of sofic shifts is of interest in symbolic dynamics.

\subsubsection{Entropy} \label{DEF:ent}

Another important aspect of shifts which is studied in dynamics is entropy (both topological and measure theoretic), though this theory is generally restricted to countable amenable groups, due to inconsistencies which exist when attempting to define entropy more generally. Notions of entropy do exist for the broader class of countable sofic groups \cite{bowen_sofic}, however only the definition for countable amenable groups is used in this paper. Formal treatment of topological and measure theoretic entropy for $G$-spaces, as well as results about these notions of entropy, can be found in \cite{kerr_li}.

\begin{defn} \label{DEF:top_ent}
    Let $G$ be a countable amenable group. Then the \textit{(topological) entropy} of a nonempty $G$-shift $X$ is defined as
    \[
        h(X) = \inf_n \frac{\log\bigl(|\cL_{F_n}(X)|\bigr)}{|F_n|} = \lim_{n \to \infty} \frac{\log\bigl(|\cL_{F_n}(X)|\bigr)}{|F_n|},
    \]
    where $\{F_n\}$ is some F{\o}lner sequence for $G$. This limit always exists and is equal to the infimum above \cite[Section 9.9]{kerr_li}. The entropy of $X$ is also independent of the choice of F{\o}lner sequence.
\end{defn}

Furthermore, some results pertain to the set of real numbers which are attained as the (topological) entropies for SFTs on a particular group.

\begin{defn}
    Let $G$ be a countable amenable group. Then let
    \[
        \cE(G) = \bigl\{h(X) : X \text{ a nonempty } G\text{-SFT}\bigr\}
    \]
\end{defn}

Note that $\cE(G)$ is a countable subset of $[0, \infty)$, since there are only countably many $G$-SFTs for any group $G$. Determining exactly what the set $\cE(G)$ is for a given group $G$ is in general quite difficult. A classic result of Lind \cite{lind_Z_entropies} precisely classifies $\cE(\Z)$ as non-negative rational multiples of logarithms of Perron numbers. More recently, Hochman and Meyerovich determined that $\cE(\Z^d)$ is the set of non-negative upper semi-computable real numbers \cite{hochman_meyerovitch}. For finitely generated amenable groups $G$ with decidable word problem which admit a translation like action by $\Z^2$, recent work by Barbieri \cite{barbieri} has classified $\cE(G)$ as the set of non-negative upper semi-computable real numbers.

With entropy, we may also define the following notion of minimality.

\begin{defn} \label{DEF:ent_min}
    Let $G$ be a countable amenable group, and $X$ a $G$-shift. Then $X$ is \textit{entropy minimal} if for each subshift $Y \subsetneq X$, we have $h(Y) < h(X)$.
\end{defn}

A weaker but related notion of minimality is SFT-entropy minimality.

\begin{defn}
    Let $G$ be a countable amenable group, and $X$ a $G$-shift. Then $X$ is \textit{SFT-entropy minimal} if for each SFT $Y \subsetneq X$, we have $h(Y) < h(X)$.
\end{defn}

Although in general SFT-entropy minimality is weaker than entropy minimality, they are in fact equivalent if the shift in question is an SFT. Proving this is a fairly standard argument involving approximating subshifts by SFTs, so we omit its proof. This fact is quite useful for proving that an SFT is entropy minimal, as it significantly reduces the amount of shifts to consider when proving entropy minimality. 

Along with topological entropy, measure-theoretic entropy can be defined if the shift $X$ is additionally endowed with a Borel probability measure (that is always Radon, since $\cA^G$ is metrizable when $G$ is countable) that behaves nicely with the shift action of $G$.

\begin{defn}
    Let $G$ be a countable amenable group, and let $X$ be a $G$-shift. A measure $\mu$ on $X$ is \textit{$G$-invariant} if for any $g \in G$ and measurable $E \subset X$, it is the case that $\mu(\sigma^{g^{-1}}E) = \mu(E)$.
    
    Let $\cM(X)$ denote the set of all $G$-invariant Borel probability measures $\mu$ on $X$.
\end{defn}

For a $G$-shift $X$ and $w \in \cL(X)$, $\mu[w]$ is used as a shorthand for $\mu([w] \cap X)$. To define the $\mu$-entropy of $X$, first an associated partition entropy must be defined.

\begin{defn}
    Let $G$ be a countable amenable group, $X$ be a $G$-shift, and $\mu \in \cM(X)$. Then, for any $F \Subset G$, the \textit{$(F,\mu)$-entropy} of $X$ is defined as
    \[
        H_\mu(X, F) = - \sum_{w \in \cL_F(X)} \mu[w] \log( \mu[w] ),
    \]
    where $0 \cdot \log(0)$ is taken to be $0$ by convention. The maximum of $H_\mu(X, F)$ over $\cM(X)$ is $\log(|\cL_F(X)|)$, and is attained only by any $\mu \in \cM(X)$ for which $\mu[w] = \frac{1}{|\cL_F(X)|}$ for all $w \in \cL_F(X)$ \cite[Corollary 4.2.1]{walters}.
\end{defn}

With this, the measure theoretic entropy can be defined.

\begin{defn} \label{DEF:mes_ent}
    Let $G$ be a countable amenable group, $X$ be a $G$-shift, and $\mu \in \cM(X)$. Then for any F{\o}lner sequence $\{F_n\}_{n=1}^\infty$ for $G$, the \textit{$\mu$-entropy} of $X$ is defined as
    \[
        h_\mu(X) = \inf_n \frac{H_\mu(X, F_n)}{|F_n|} = \lim_{n \to \infty} \frac{H_\mu(X, F_n)}{|F_n|}.
    \]
    As with topological entropy, this limit always exists, is equal to this infimum, and is independent of the choice of F{\o}ner sequence \cite[Section 9.3]{kerr_li}.
    Furthermore, the Variational Principle \cite[Theorem 9.43]{kerr_li} gives that
    \[
        h(X) = \sup_{\mu \in \cM(X)} h_\mu(X).
    \]
    A measure $\mu \in \cM(X)$ satisfying $h(X) = h_\mu(X)$ is called a \textit{measure of maximal entropy}, and for $G$-shifts, there always exists at least one measure of maximal entropy, since shift actions are expansive.
\end{defn}

\section{Free Extension Shifts} \label{FE}

Though the primary purpose of this paper is to prove that locally finite groups are precisely the groups which exhibit strong dynamical properties for all SFTs, proving many of these properties directly is somewhat tedious. Instead, we develop a general theory of \textit{free extension} shifts, which simplifies (and even trivializes) many of the results for locally finite groups. Essentially all of the primary results in this paper use properties of free extensions, which are constructed in this section.

The notion of a free extension shift is not new however. Hochman and Meyerovich \cite{hochman_meyerovitch} used them (though not explicitly by name) in their landmark paper characterizing the possible entropies of $\Z^d$ SFTs. The term free extension and some associated notation used were coined by Barbieri \cite{barbieri}, with free extensions appearing as a special case of a far more general method of constructing \dquote{extensions} of shifts. The construction given here is far less general than Barbieri's construction, but it is perhaps more amenable to specifically analyzing free extensions.

\subsection{Construction of Free Extensions}

We construct free extension shifts by taking a homeomorphic image of particular closed, $G$-invariant subset of a $G$-space which is conjugate to the full $G$-shift. This $G$-space is not a $G$-shift in the typical sense however, as the action of $G$ on this space is not given by the typical shifts. By explicitly constructing a conjugacy between these two $G$-spaces, we can easily construct subsystems of this $G$-space and push them forward under this conjugacy into a $G$-shift.

To begin, we need to first define the notion of choice functions for the collection of right cosets of a subgroup.

\begin{defn}
    Let $G$ be a group and $H \le G$. A \textit{choice function} for $H \backslash G$ (the right cosets of $H$ in $G$) is any function $\zeta:H \backslash G \to G$ such that for any $Hc \in H \backslash G$, $\zeta(Hc) \in Hc$. Since the cosets of a subgroup are mutually disjoint, any choice function is injective. Restricting the codomain of a choice function $\zeta$ to $\zeta(H \backslash G)$ results in each choice function being a bijection. For simplicity, we also let $C_\zeta = \zeta(H \backslash G)$ so that $\zeta:H \backslash G \to C_\zeta$ is a bijection.
    
    We take $\cC(H \backslash G)$ to be the set of all choice functions $\zeta:H \backslash G \to C_\zeta$, so that each are bijections. Note that by the axiom of choice, choice functions always exist, regardless of the choice of $G$ and $H$.
\end{defn}

With these choice functions for a group $G$ and subgroup $H$, we can now define the $G$-space which is conjugate to the full $G$-shift (see Theorem \ref{FE:cons_conj}).

\begin{defn}
    Let $G$ be a group, $H \le G$, and $\cA$ be a finite alphabet. Let $\zeta \in \cC(H \backslash G)$. $G$ admits a natural action on the collection of right cosets $H \backslash G$ by $\eta^g(Ha) = Hag$. Note that each $\eta^g$ is a well defined bijection from $H \backslash G$ to itself. With these two maps, we may define an action of $G$ on $C_\zeta$, where for each $g \in G$, the action of $g$ on $C_\zeta$ is $\xi_{\zeta}^g:C_\zeta \to C_\zeta$, which is defined as $\xi_{\zeta}^g = \zeta \circ \eta^g \circ \zeta^{-1}$. This makes $(H \backslash G, \eta)$ and $(C_\zeta, \xi_\zeta)$ conjugate $G$-spaces when $H \backslash G$ and $C_\zeta$ are endowed with the discrete topology.
    
    Let $\cA^{H ; C_\zeta}$ denote the set $\bigl(\cA^H\bigr)^{C_\zeta}$. This space is endowed with the product topology from the product topology on $\cA^H$, which makes $\cA^{H ; C_\zeta}$ a compact Hausdorff space (which is metrizable when $G$ is countable). For each $g \in G$, define $\rho_\zeta^g:\cA^{H ; C_\zeta} \to \cA^{H ; C_\zeta}$ by
    \[
        \rho_\zeta^g\bigl(\{w_c\}_{c \in C_\zeta}\bigr) = \bigl\{\sigma^{cg\xi_{\zeta}^g(c)^{-1}} w_{\xi_{\zeta}^g(c)}\bigr\}_{c \in C_\zeta}.
    \]
    Since $\xi_{\zeta}^g \circ \zeta = \zeta \circ \eta^g$, it follows that for any $c \in C_\zeta$,
    \[
        \xi_{\zeta}^g(c) = \xi_{\zeta}^g(\zeta(Hc)) = \zeta(\eta^g(Hc)) = \zeta(Hcg),
    \]
    and therefore $Hcg = H\xi_{\zeta}^g(c)$. This gives that $cg\xi_{\zeta}^g(c)^{-1} \in H$, and so $\rho_\zeta^g$ is well defined. In fact, each $\rho_\zeta^g$ is a homeomorphism because the reordering of the indices only depends on $g$, and each $\sigma$ is a homeomorphism on $\cA^H$. The collection $\rho_\zeta = \{\rho_\zeta^g : g \in G\}$ forms an action of $G$ on $\cA^{H ; C_\zeta}$, which gives that $(\cA^{H ; C_\zeta}, \rho_\zeta)$ is a $G$-space.
\end{defn}

Informally, the $G$-space $(\cA^{H ; C_\zeta}, \rho_\zeta)$ can be thought of as having copies of $\cA^H$ on each coset of $H$ in $G$, with elements of $\cA^{H ; C_\zeta}$ having a separate configuration from $\cA^H$ on each coset. Dynamically, $\rho_\zeta$ acts on this space in two steps. The first is that for a particular choice of $g \in G$, the configurations in each coset get shuffled to other cosets (although it is certain that one configuration remains in the coset it originated in), based on how $g$ acts on the cosets of $H$ in $G$. The second step is to then individually shift each configuration in its new coset by an element of $H$, which is determined from the coset the configuration is now in, and the particular choice of $g$. The dynamics are completely independent of the particular element of $\cA^{H ; C_\zeta}$. The particular choice for the definition of $\rho_\zeta$ results in the dynamics of $(\cA^{H ; C_\zeta}, \rho_\zeta)$ being conjugate to the full $G$-shift under the $\sigma$-action of $G$ for any choice of $\zeta$ (see Theorem \ref{FE:cons_conj} below). Other similar systems can be constructed on the set $\cA^{H; B}$ for any set $B$ by choosing functions $\alpha:G \times B \to B$ and $\beta:G \times B \to H$ with $\alpha(g, \cdot)$ a bijection for each $g \in G$, and defining an action $\delta_{\alpha,\beta}$ of $G$ on $\cA^{H ; B}$ by
\[
    \delta^g_{\alpha,\beta}\bigl(\{w_b\}_{b \in B}\bigr) = \bigl\{\sigma^{\beta(g,b)}w_{\alpha(g,b)}\}_{b \in B}.
\]
Such systems may be interesting to study in their own right, however we only focus on the particular systems $(\cA^{H; C_\zeta}, \rho_\zeta)$. Let us now define the conjugacy map between these spaces.

\begin{defn}
    Let $G$ be a group, let $H \le G$, let $\zeta \in \cC(H \backslash G)$, and let $\cA$ be a finite alphabet. Define a \textit{construction} function $\kappa_\zeta:\cA^{H ; C_\zeta} \to \cA^G$, by
    \[
        \kappa_\zeta(\{w_c\}_{c \in C_\zeta}) = \bigvee_{c \in C_\zeta} \sigma^{c^{-1}}w_c.
    \]
    Note that $\kappa_\zeta$ is well defined since $\sigma^{c^{-1}}w_c$ has shape $Hc$, so each pattern in the join is disjoint, and  furthermore, $\kappa_\zeta$ produces a configuration in $\cA^G$ since $\{Hc\}_{c \in C_\zeta}$ is a partition of $G$.
\end{defn}

It remains to show that any choice of construction function $\kappa_\zeta$ is in fact a conjugacy between $(\cA^{H; C_\zeta}, \rho_\zeta)$ and $(\cA^G, \sigma)$.

\begin{thm} \label{FE:cons_conj}
    Let $G$ be group, let $H \le G$, and let $\zeta \in \cC(H \backslash G)$. Then $\kappa_\zeta$ is a conjugacy between $(\cA^{H ; C_\zeta}, \rho)$ and $(\cA^G, \sigma)$. That is, $\kappa_\zeta$ is a homeomorphism, and $\forall g \in G$,
    \[
        \kappa_\zeta\bigl(\rho_\zeta^g(\{w_c\}_{c \in C_\zeta})\bigr) = \sigma^g\bigl(\kappa_\zeta(\{w_c\}_{c \in C_\zeta})\bigr).
    \]
    \begin{proof}
        First, we show that $\kappa_\zeta$ is injective. Let $\{w_c\}_{c \in C_\zeta}, \{v_c\}_{c \in C_\zeta} \in \cA^{H ; C_\zeta}$ be such that $\kappa_\zeta(\{w_c\}) = \kappa_\zeta(\{v_c\})$. For any $d \in C_\zeta$, it must be that $\kappa_\zeta(\{w_c\})|_{Hd} = \kappa_\zeta(\{v_c\})|_{Hd}$. As $d \in C_\zeta$, it must be that $\zeta(Hd) = d$, and so $\kappa_\zeta(\{w_c\})|_{Hd} = \sigma^{d^{-1}}w_d$, which gives $\sigma^{d^{-1}}w_d = \sigma^{d^{-1}}v_d$. This implies $w_d = v_d$, and since $d \in C_\zeta$ was arbitrary, this gives $\{w_c\}_{c \in C_\zeta} = \{v_c\}_{c \in C_\zeta}$.
        
        Next, let us show that $\kappa_\zeta$ is surjective. Let $x \in \cA^G$. Then
        \begin{align*}
            x &= \bigvee_{c \in C_\zeta} x|_{Hc} = \bigvee_{c \in C_\zeta} \sigma^{c^{-1}}\bigl(\sigma^{c}(x|_{Hc})\bigr) \\
            &= \bigvee_{c \in C_\zeta} \sigma^{c^{-1}}\bigl((\sigma^{c}x)|_H\bigr) = \kappa_\zeta\bigl(\{(\sigma^cx)|_H\}_{c \in C_\zeta}\bigr).
        \end{align*}
        
        Now we show that $\kappa_\zeta$ is continuous. It suffices to show that for any $w \in \cL(\cA^G)$, the set $\kappa_\zeta^{-1}([w]_G)$ is open in $\cA^{H ; C_\zeta}$. Let $F$ be the shape of $w$, and for each $c \in C_\zeta$, let $F_c = F \cap Hc$. Note that since $F$ is finite, at most finitely many $F_c$ are nonempty.
        
        Recall that for any sets $E \subset F$, and $w \in \cA^E$, the $F$-extensions of $w$ are defined to be $[w]_F = \{v \in \cA^F : v|_E = w\}$, and let $U = \prod_{c \in C_\zeta} \bigl[\sigma^{c}(w|_{F_c})\bigr]_H$. For any $c \in C_\zeta$ with $F_c$ empty (which occurs for all but finitely many $c \in C_\zeta$), $[\sigma^{c}(w|_{F_c})]_H = \cA^H$ by definition. As such, $U$ is a basic open set in $\cA^{H ; C_\zeta}$. We now show that $\kappa_\zeta^{-1}([w]_G) = U$. Since $F_c \subset Hc$, it follows that $x|_{Hc} \in \bigl[x|_{F_c}\bigr]_{Hc}$ for any $x \in \cA^G$. Using this, for any $x \in [w]_G$ we have $x|_{Hc} \in \bigl[x|_{F_c}\bigr]_{Hc} = \bigl[w|_{F_c}\bigr]_{Hc}$, and for any $c \in C_\zeta$,
        \[
            \bigl(\kappa_\zeta^{-1}(x)\bigr)_c = (\sigma^cx)|_H = \sigma^c(x|_{Hc}) \in \sigma^{c}[w|_{F_c}]_{Hc} = \bigl[\sigma^{c}(w|_{F_c})\bigr]_{H},
        \]
        and therefore $\kappa_\zeta^{-1}(x) \in U$. Now let $y = \{y_c\}_{c \in C_\zeta} \in U$. With $x = \kappa_\zeta(y)$, it is the case that $x \in [w]_G$ if and only if $x|_{F_c} = w|_{F_c}$ for all $c \in C_\zeta$. Indeed, for each $c \in C_\zeta$, 
        \[
            x|_{Hc} = \sigma^{c^{-1}}y_c \;\in\; \sigma^{c^{-1}}\bigl[\sigma^{c}(w|_{F_c})\bigr]_H = \bigl[\sigma^{c^{-1}c}(w|_{F_c})\bigr]_{Hc} = [w|_{F_c}]_{Hc},
        \]
        and therefore $x|_{F_c} = w|_{F_c}$, which gives $x \in [w]_G$. This gives that $\kappa_\zeta(U) \subset [w]_G$, and therefore $\kappa_\zeta^{-1}([w]_G) = U$, so $\kappa_\zeta^{-1}([w]_G)$ is open.
        
        Using this, $\kappa_\zeta$ is a homeomorphism, because $\kappa_\zeta$ is a continuous bijection from $\cA^{H ; C_\zeta}$, which is compact, to $\cA^G$, which is Hausdorff.
        
        Finally, for any $\{w_c\}_{c \in C_\zeta} \in \cA^{H ; C_\zeta}$ and $g \in G$,
        \begin{align*}
            \kappa_\zeta\bigl(\rho_\zeta^g(\{w_c\}_{c \in C_\zeta})\bigr) 
            &= \kappa_\zeta\bigl(\{\sigma^{cg\xi_{\zeta}^g(c)^{-1}}w_{\xi_{\zeta}^g(c)}\}_{c \in C_\zeta}\bigr) \\
            &= \bigvee_{c \in C_\zeta} \sigma^{c^{-1}}\Bigl(\sigma^{cg\xi_{\zeta}^g(c)^{-1}}w_{\xi_{\zeta}^g(c)}\Bigr) \\
            &= \sigma^g \left(\bigvee_{c \in C_\zeta} \sigma^{\xi_{\zeta}^g(c)^{-1}}w_{\xi_{\zeta}^g(c)}\right) \\
            &\stackrel{(a)}{=} \sigma^g \left(\bigvee_{c \in C_\zeta} \sigma^{c^{-1}} w_c\right) \\
            &= \sigma^g\bigl(\kappa_\zeta(\{w_c\}_{c \in C_\zeta})\bigr),
        \end{align*}
        where equality $(a)$ follows from reindexing the join, which is possible because $\xi_{\zeta}^g$ is a permutation of $C_\zeta$.
    \end{proof}
\end{thm}

Although this conjugacy is well defined, it is entirely dependent on the choice of $\zeta$. Due to the correspondence between $C_\zeta$ and $H \backslash G$ however, while the precise definition of the conjugacy depends on $\zeta$, the following result gives that the objects we define with the construction function are not dependent on $\zeta$. This makes it possible to define free extensions independently of $\zeta$, which is done in Definition \ref{FE:defn}.

\begin{lemma} \label{FE:unique_ext}
    Let $G$ be a group, let $H \le G$, and let $\zeta_1, \zeta_2 \in \cC(H \backslash G)$. Let $\kappa_1 = \kappa_{\zeta_1}$ and $\kappa_2 = \kappa_{\zeta_2}$, and let $C_1 = C_{\zeta_1}$ and $C_2 = C_{\zeta_2}$. Then for any $H$-shift $X$, $\kappa_1(X^{C_1}) = \kappa_2(X^{C_2})$.
    \begin{proof}
        It suffices to show that $\kappa_1(X^{C_1}) \subset \kappa_2(X^{C_2})$, since the roles of $\zeta_1$ and $\zeta_2$ may be reversed to deduce the reverse inclusion. 
        
        Let $x \in \kappa_1(X^{C_1})$, and let $\{w_c\}_{c \in C_1} = \kappa_1^{-1}(x)$, which is well defined because $x$ is in the image of $\kappa_1$, and $\kappa_1$ is injective by Theorem \ref{FE:cons_conj}. Now, let $f = \zeta_1 \circ \zeta_2^{-1}$, which is a bijection from $C_2$ to $C_1$, and note that for any $c \in C_2$, we have that $f(c) \in Hc$, which implies $cf(c)^{-1} \in H$. Since for each $c \in C_2$ we have $w_{f(c)} \in X$, and $X$ is an $H$-shift, it follows that $\sigma^{cf(c)^{-1}}w_{f(c)} \in X$. Then
        \begin{align*}
            \kappa_2(\{\sigma^{cf(c)^{-1}}w_{f(c)}\}_{c \in C_2}) &= \bigvee_{c \in C_2} \sigma^{c^{-1}}\bigl(\sigma^{cf(c)^{-1}}w_
            {f(c)}\bigr) = \bigvee_{c \in C_2} \sigma^{f(c)^{-1}} w_{f(c)} \\
            &\stackrel{(a)}{=} \bigvee_{c \in C_1} \sigma^{c^{-1}} w_c = \kappa_1(\{w_c\}) = x,
        \end{align*}
        where equality $(a)$ follows from reindexing the join, which is possible since $f$ is a bijection from $C_2$ to $C_1$.
        Therefore, $\kappa_1(X^{C_1}) \subset \kappa_2(X^{C_2})$.
    \end{proof}
\end{lemma}

With this result, we may finally formally define free extension shifts.

\begin{defn} \label{FE:defn}
    Let $G$ be a group, let $H \le G$, and let $X$ be an $H$-shift, called the \textit{base} shift. Let $\zeta \in \cC(H \backslash G)$. The \textit{free $G$-extension of $X$}, denoted $X^{\uparrow G}$, is defined to be $\kappa_\zeta(X^{C_\zeta})$, which is independent of the choice of $\zeta$ by Lemma \ref{FE:unique_ext}. $X^{\uparrow G}$ is a $G$-shift, because $X^{C_\zeta}$ is $\rho_\zeta$-invariant due to $X$ being $H$-invariant, and closed because $X$ is closed.
    
    Since the free extension of $X$ is independent of $\zeta$, let $\cA^{H : G} = \cA^{H ; H \backslash G} = \bigl(\cA^H\bigr)^{H \backslash G}$, which is homeomorphic to $\cA^{H ; C_\zeta}$ for any choice of $\zeta$, since $(C_\zeta, \xi_\zeta)$ is conjugate to $(H \backslash G, \eta)$. We then consider a generic construction function $\kappa:\cA^{H : G} \to \cA^G$ as being $\kappa_\zeta$ for some arbitrary $\zeta \in \cC(H \backslash G)$. Using this notation, $X^{\uparrow G}$ may also be written as $X^{\uparrow G} = \kappa(X^{H \backslash G})$. Similarly, we define $C_\kappa$ to be $C_\zeta$ for the arbitrary choice of $\kappa = \kappa_\zeta$, $\rho_\kappa$ to be $\rho_\zeta$, and $\xi^g_\kappa = \xi^g_\zeta$.
    
    In the case that $H = G$, it is trivial to see that for any $H$-shift $X$, we have $X^{\uparrow G} = X$. Most of the properties in the next section are trivial in this case, and arguments are made primarily in the case that $H$ is a proper subgroup of $G$.
\end{defn}

\subsection{Properties of Free Extensions and their Base Shifts}

Now that we have established that free extensions are well defined, we now prove that many useful properties can be transferred from a base shift to its extension, and vice-versa.

First, we observe that the topological entropy of a free extension is the same as the entropy of the shift being extended.

\begin{prop}[Proposition 5.2, \cite{barbieri}] \label{FE:entropy}
    Let $G$ be a countable amenable group, $H \le G$, and $X$ an $H$-shift. Then
    \[
        h(X^{\uparrow G}) = h(X).
    \]
\end{prop}

Second, if there are 3 groups $K \le H \le G$, then taking a $K$-shift and extending it to $G$ produces the same thing as first extending to $H$, and then to $G$.

\begin{lemma} \label{FE:mult_ext}
    Let $G$ be a group, $K \le H \le G$, and $Y$ a $K$-shift. Then $Y^{\uparrow G} = (Y^{\uparrow H})^{\uparrow G}$.
    \begin{proof}
        Let $\zeta_1 \in \cC(K \backslash H)$ and let $\zeta_2 \in \cC(H \backslash G)$, and let $C_1 = C_{\zeta_1}$ and $C_2 = C_{\zeta_2}$. Note that $C_3 = C_1C_2 \subset G$ contains exactly one element from each coset in $K \backslash G$. This means we can identify elements of $C_1 \times C_2$ with $C_3$, with the natural bijection $(c,d) \mapsto cd$. As such, we may define $\zeta_3 \in \cC(K \backslash G)$ such that $C_{\zeta_3} = C_3$. Let $\kappa_i = \kappa_{\zeta_i}$ for $i = 1, 2, 3$. Then for any $\{w_c\}_{c \in C_3} \in Y^{C_3}$,
        \begin{align*}
            \kappa_3(\{w_c\}) &= \bigvee_{c \in C_3} \sigma^{c^{-1}}w_{c} \\
            &= \bigvee_{(c,d) \in C_1 \times C_2} \sigma^{(cd)^{-1}}w_{c,d} \\
            &= \bigvee_{d \in C_2} \sigma^{d^{-1}} \bigvee_{c \in C_1} \sigma^{c^{-1}} w_{c,d} \\
            &= \bigvee_{d \in C_2} \sigma^{d^{-1}}\kappa_1(\{w_{c,d}\}_{c \in C_1}) \\
            &= \kappa_2\Bigl(\bigl\{\kappa_1(\{w_{c,d}\}_{c \in C_1})\bigr\}_{d \in C_2}\Bigr),
        \end{align*}
        and therefore $Y^{\uparrow G} \subset \bigl(Y^{\uparrow H}\bigr)^{\uparrow G}$ by Lemma \ref{FE:unique_ext}. Similarly, starting with some arbitrary $\bigl\{\{w_{c,d}\}_{c \in C_1}\bigr\}_{d \in C_2} \in (Y^{K \backslash H})^{H \backslash G}$, and using the equality above in reverse, it follows that $\bigl(Y^{\uparrow H}\bigr)^{\uparrow G} \subset Y^{\uparrow G}$, and therefore they are equal.
    \end{proof}
\end{lemma}

Next, we prove a stability result for free extensions, namely that the intersection of free extensions is the free extension of an intersection of shifts. This is also the first proof where the particular choice of $\zeta$ is omitted, and an arbitrary construction function is used.

\begin{lemma} \label{FE:intersect}
    Let $G$ be a group, $H \le G$, and $\{Y_i\}_{i \in I}$ be a collection of $H$-shifts. Then $\bigcap_{i \in I} Y_i^{\uparrow G} = \left(\bigcap_{i \in I} Y_i\right)^{\uparrow G}$.
    \begin{proof}
        Let $\kappa$ be an arbitrary construction function on $(\cA^{H : G}, \rho_\kappa)$. We clearly have that
        \[
            \bigcap_{i \in I} Y_i^{H : G} = \left(\bigcap_{i \in I} Y_i\right)^{H : G},
        \]
        and that both of these are closed and $\rho_\kappa$-invariant. $\kappa$ is a bijective homeomorphism by Theorem \ref{FE:cons_conj}, so we have
        \[
            \kappa\left(\bigcap_{i \in I} Y_i^{H : G}\right) = \kappa\left( \left(\bigcap_{i \in I} Y_i\right)^{H : G}\right) = \left(\bigcap_{i \in I} Y_i\right)^{\uparrow G}.
        \]
        Again using that $\kappa$ is a bijection, we have that
        \[
            \kappa\left(\bigcap_{i \in I} Y_i^{H : G}\right) = \bigcap_{i \in I} \kappa(Y_i^{H : G}) = \bigcap_{i \in I} Y_i^{\uparrow G}.
        \]
        Combining the results of the two previous displays, we obtain the desired result.
    \end{proof}
\end{lemma}

The fourth result is fairly natural and deals with taking free extensions defined by a set of forbidden patterns. Essentially, for a selection of forbidden patterns on a group $H \le G$, taking the $H$-shift defined by the forbidden patterns and then extending this shift to $G$ results in the same shift as the $G$-shift defined by the forbidden patterns directly.

\begin{lemma} \label{FE:forbidden_ext}
    Let $G$ be a group and $H \le G$. Then for any $\bF \subset \cL(\cA^H)$, $\cX^G[\bF] = (\cX^H[\bF])^{\uparrow G}$.
    \begin{proof}
        Let $\kappa$ be a construction function on $\cA^{H : G}$, and let $d \in C_\kappa$ such that $H = Hd$. Note that for any $w \in \bF$ with shape $F$, we have $F \Subset Hd$. Let $x \in (\cX^H[\bF])^{\uparrow G}$, and $\{w_c\}_{c \in C_\kappa} = \kappa^{-1}(x)$. Then for any $g \in G$ and $F \Subset H$,
        \begin{align*}
            (\sigma^gx)|_F &= \bigl(\sigma^g\kappa(\{w_c\})\bigr)|_F \\
            &= \kappa\bigl(\rho_\kappa^g(\{w_c\})\bigr)|_F \\
            &= \kappa\bigl(\{\sigma^{cg\xi_{\kappa}^g(c)^{-1}} w_{\xi_{\kappa}^g(c)}\}\bigr)|_F \\
            &= (\sigma^{d^{-1}}\sigma^{dg\xi_{\kappa}^g(d)^{-1}}w_{\xi_{\kappa}^g(d)})|_F \\
            &= (\sigma^{g\xi_{\kappa}^g(d)^{-1}}w_{\xi_{\kappa}^g(d)})|_F.
        \end{align*}
        Since $w_{\xi_{\kappa}^g(d)} \in \cX^H[\bF]$, so is $\sigma^{g\xi_{\kappa}^g(d)^{-1}} w_{\xi_{\kappa}^g(d)}$ (since $d^{-1} \in H$ and $dg\xi_{\kappa}^g(d)^{-1} \in H$, we have $g\xi_{\kappa}^g(d)^{-1} \in H$), and therefore it must be that $(\sigma^gx)|_F \notin \bF$. As such, $(\cX^H[\bF])^{\uparrow G} \subset \cX^G[\bF]$.
        
        Now, for any $x \in \cX^G[\bF]$, let $\{w_c\}_{c \in C_\kappa} = \kappa^{-1}(x)$. Then for all $c \in C_\kappa$, we have $w_c = (\sigma^cx)|_H$, so for every $h \in H$ and $F \Subset H$,
        \[
            \bigl(\sigma^h(\sigma^cx)|_H\bigr)|_F = \sigma^h\Bigl(\bigl((\sigma^cx)|_H\bigr)|_{Fh}\Bigr) = \sigma^h\bigl((\sigma^cx)|_{Fh}\bigr) = (\sigma^{hc}x)|_F \notin \bF,
        \]
        therefore $w_c \in \cX^H[\bF]$. As such, $\cX^G[\bF] \subset (\cX^H[\bF])^{\uparrow G}$, which, when combined with the result of the previous paragraph, gives $(\cX^H[\bF])^{\uparrow G} = \cX^G[\bF]$.
    \end{proof}
\end{lemma}

Using this result, it is easy to see that the free extension of an SFT remains an SFT. Perhaps unsurprisingly, the converse also holds; if the free extension of a shift is an SFT, then the base shift must have been an SFT to begin with.

\begin{lemma} \label{FE:sft}
    Let $G$ be a group, $H \le G$, and $Y$ be an $H$-shift. Then $Y$ is an SFT if and only if $Y^{\uparrow G}$ is an SFT.
    \begin{proof}
        First, suppose $Y$ is an SFT. Let $F \Subset H$ be a forbidden shape for $Y$, and $\bF \Subset \cA^F$ be a set of forbidden $F$-patterns so that $Y = \cX^H[\bF]$. By Lemma \ref{FE:forbidden_ext}, $Y^{\uparrow G} = (\cX^H[\bF])^{\uparrow G} = \cX^G[\bF]$, which is clearly an SFT.
        
        Now suppose that $Y^{\uparrow G}$ is an SFT. Let $\kappa$ be a construction function on $\cA^{H : G}$, and let $d \in C_\kappa$ such that $H = Hd$. Let $F \Subset G$ be a forbidden shape for $Y^{\uparrow G}$ so that $Y^{\uparrow G} = \cX^G[\cF_F(Y^{\uparrow G})]$. Now, define $C_0 \subset C_\kappa$ to be the set of all $c \in C_\kappa$ for which $F \cap Hc \ne \varnothing$, and for $c \in C_0$, let $F_c = F \cap Hc$. This partitions $F$ into the finitely many disjoint subsets $F_c$, which are each contained within a separate coset of $H$. Then define
        \[
            E = \bigcup_{c \in C_0} F_cc^{-1} \subset H,
        \]
        and
        \[
            \hat{F} = \bigcup_{c \in C_0} Ec.
        \]
        Note that for each $c \in C_0$ we have $F_cc^{-1} \subset E$, and therefore $F_c = (F_cc^{-1})c \subset Ec \subset \hat{F}$, so $F \subset \hat{F}$. As such, $Y^{\uparrow G} = \cX[\cF_{\hat{F}}(Y^{\uparrow G})]$. For any $w \in \cA^E$, let
        \[
            P(w) = \bigcup_{c \in C_0} [\sigma^{c^{-1}}w]_{\hat{F}} \subset \cA^{\hat{F}},
        \]
        and define
        \[
            \bF = \bigl\{w \in \cA^E : P(w) \subset \cF_{\hat{F}}(Y^{\uparrow G})\bigr\}.
        \]
        Let us now show that $Y = \cX^H[\bF]$, which clearly shows $Y$ is a $H$-SFT, since $\bF$ is finite.
        
        First, let $x \in Y^{\uparrow G} = \cX^G[\cF_{\hat{F}}(Y^{\uparrow G})]$. Let $g \in G$, and pick any $c \in C_0$. Then, from the definition of $P(w)$,
        \[
            (\sigma^{c^{-1}g}x)|_{\hat{F}} \in [(\sigma^{c^{-1}g}x)|_{Ec}]_{\hat{F}} =  [\sigma^{c^{-1}}(\sigma^gx)|_E]_{\hat{F}} \subset P\bigl((\sigma^gx)|_E\bigr).
        \]
        Since $x \in \cX^G[\cF_{\hat{F}}(Y^{\uparrow G})]$, it follows that $(\sigma^{c^{-1}g}x)|_{\hat{F}} \notin \cF_{\hat{F}}(Y^{\uparrow G})$, and therefore, $P\bigl((\sigma^gx)|_E\bigr) \not\subset \cF_{\hat{F}}(Y^{\uparrow G})$. This gives that $(\sigma^{g}x)|_E \notin \bF$. Since $g$ was arbitrary, this implies $x \in \cX^G[\bF]$, and therefore $Y^{\uparrow G} \subset \cX^G[\bF]$.
        
        Now let $x \in \cX^G[\bF]$. By definition, for each $g \in G$, it must be that $(\sigma^gx)|_E \notin \bF$, and therefore for every $c \in C_0$, we have $P\bigl((\sigma^{cg}x)|_E\bigr) \not\subset \cF_{\hat{F}}(Y^{\uparrow G})$. As such, for each $c \in C_0$, there exists $w_c \in P\bigl((\sigma^{cg}x)|_E\bigr) \setminus \cF_{\hat{F}}(Y^{\uparrow G})$. Note that this means for each $c \in C_0$, there exists $d_c \in C_0$ such that $(\sigma^{d_{c}}w_c)|_E = (\sigma^{cg}x)|_E$. Let $x_c \in Y^{\uparrow G}$ be such that $x_c|_{\hat{F}} = w_c$. This must be possible since $w_c \notin \cF_{\hat{F}}(Y^{\uparrow G})$, and $Y^{\uparrow G}$ is shift invariant. Let $\{a_{c,d}\}_{d \in C_\kappa} = \kappa^{-1}(x_c)$, and note that $a_{c,d} \in Y$ for each $c \in C_0$ and $d \in C_\kappa$. Furthermore, for each $c \in C_0$, we have $a_{c, d_c}|_E = (\sigma^{cg}x)|_E$, which gives that
        \[
            (\sigma^{c^{-1}}a_{c, d_c})|_{Ec} = \sigma^{c^{-1}}(a_{c,d_c}|_E) = \sigma^{c^{-1}}\bigl((\sigma^{cg}x)|_E\bigr) = (\sigma^gx)|_{Ec}.
        \]
        Define $\{y_d\}_{d \in C_\kappa} \in Y^{H \backslash G}$ as follows. For each $c \in C_0$, let $y_c = a_{c, d_c}$, and for $d \in C \setminus C_0$, let $y_d \in Y$ (it does not matter how these are chosen). Then $y = \kappa(\{y_d\}) \in Y^{\uparrow G}$, so for each $c \in C_0$, it is the case that $y|_{Ec} = (\sigma^{c^{-1}}a_{c, d_c})|_{Ec} = (\sigma^gx)|_{Ec}$, which gives that $y|_{\hat{F}} = (\sigma^gx)|_{\hat{F}}$. Since $y \in Y^{\uparrow G}$, this implies $(\sigma^gx)|_{\hat{F}} \notin \cF_{\hat{F}}(Y^{\uparrow G})$. As this is true for all $g \in G$, this implies $x \in \cX[\cF_{\hat{F}}(Y^{\uparrow G})] = Y^{\uparrow G}$, and therefore $\cX^G[\bF] \subset Y^{\uparrow G}$.
        
        The results of the two previous paragraphs gives that $Y^{\uparrow G} = \cX^G[\bF]$. By Lemma \ref{FE:forbidden_ext}, since $\bF \subset \cA^E$, and $E \subset H$, we have $Y^{\uparrow G} = \cX^G[\bF] = (\cX^H[\bF])^{\uparrow G}$, and so $\kappa(Y^{H \backslash G}) = \kappa\bigl((\cX^H[\bF])^{H \backslash G}\bigr)$. Since $\kappa$ is a bijection, it must be that $Y^{H \backslash G} = (\cX^H[\bF])^{H \backslash G}$, and therefore $\cX^H[\bF] = Y$.
    \end{proof}
\end{lemma}

Next, we show a similar result to the previous one, replacing the property of being a $G$-SFT with being strongly irreducible.

\begin{lemma} \label{FE:si}
    Let $G$ be a group, $H \le G$, and $Y$ be an $H$-shift. Then $Y$ is strongly irreducible if and only if $Y^{\uparrow G}$ is strongly irreducible.
    \begin{proof}
        First, suppose that $Y$ is strongly irreducible. Then there exists $K \Subset H$ such that for any $u,v \in \cL(Y)$ with shapes $F_u$ and $F_v$ such that $F_u \cap KF_v = \varnothing$, then there exists $x \in Y^{\uparrow G}$ such that $x|_{F_u} = u$ and $x|_{F_v} = v$. We will now show that $X = Y^{\uparrow G}$ is strongly irreducible with the same $K$. Let $u, v \in \cL(X)$ with shapes $F_u$ and $F_v$ such that $F_u \cap KF_v = \varnothing$. Since $u,v \in \cL(X)$, let $x_u, x_v \in X$ such that $x_u|_{F_u} = u$ and $x_v|_{F_v} = v$. Let $\kappa$ be some construction function on $\cA^{H : G}$. Since $x_u, x_v \in X = Y^{\uparrow G}$, we may take $\{y_c\}_{c \in C_\kappa} = \kappa^{-1}(x_u)$ and $\{z_c\}_{c \in C_\kappa} = \kappa^{-1}(x_v)$ with $y_c, z_c \in Y$ for all $c \in C_\kappa$. Let $c \in C_{\kappa}$, and define $U_c = F_uc^{-1} \cap H$ and $V_c = F_vc^{-1} \cap H$. With $U_c$ and $V_c$ being subsets of $H$, we have that $y_c|_{U_c} \in \cL(Y)$ and $z_c|_{V_c} \in \cL(Y)$. Furthermore, we have
        \[
            U_c \cap KV_c \subset F_uc^{-1} \cap KF_vc^{-1} = (F_u \cap KF_v)c^{-1} = \varnothing.
        \]
        As such, the strong irreducibility of $Y$ gives that there exists $w_c \in Y$ such that $w_c|_{U_c} = y_c|_{U_c}$ and $w_c|_{V_c} = z_c|_{V_c}$. Let $x = \kappa(\{w_c\}_{c \in C_\kappa}) \in X$. We now show that $x|_{F_u} = u$ and $x|_{F_v} = v$. Indeed, as
        \[
            F_u = \bigsqcup_{c \in C_\kappa} F_u \cap Hc = \bigsqcup_{c \in C_\kappa} (F_uc^{-1} \cap H)c = \bigsqcup_{c \in C_\kappa} U_cc,
        \]
        and similarly for $F_v$ with $V_c$ in place of $U_c$, we may obtain that $x|_{F_u} = u$ by checking $x|_{U_cc} = u|_{U_cc}$ for every $c \in C_\kappa$ (and similarly for $x|_{F_v} = v$). For $c \in C_\kappa$, we have $U_cc \subset Hc$, and so
        \begin{align*}
            x|_{U_cc} = (x|_{Hc})|_{U_cc} &= (\sigma^{c^{-1}}x_c)|_{U_cc} = \sigma^{c^{-1}}(x_c|_{U_c}) = \sigma^{c^{-1}}(y_c|_{U_c}) \\
            &= \sigma^{c^{-1}}\bigl((\sigma^cx_u)|_{U_c}\bigr) = \sigma^{c^{-1}}\sigma^c(x_u|_{U_cc}) = u|_{U_cc},
        \end{align*}
        where the final equality follows from the fact that $x_u|_{F_u} = u$. By the same argument, replacing $U_c$ with $V_c$, $y_c$ with $z_c$, $x_u$ with $x_v$, and $u$ with $v$, we obtain $x|_{V_cc} = v|_{V_cc}$, and therefore $x|_{F_u} = u$ and $x|_{F_v} = v$. As such, $Y^{\uparrow G}$ is strongly irreducible.

        Now, suppose that $X = Y^{\uparrow G}$ is strongly irreducible. Let $L \Subset G$ be such that if $u, v \in \cL(X)$ with shapes $E_u$ and $E_v$ satisfy $E_u \cap LE_v = \varnothing$, then there exists $x \in X$ such that $x|_{E_u} = u$ and $x|_{E_v} = v$. Let $K = L \cap H$, and we will now show that $Y$ is strongly irreducible with this $K$. Note that $K$ must be nonempty because $L$ must contain an element of $H$ (otherwise, if $E_u$ and $E_v$ are finite subsets of $H$ which intersect, then $E_u \cap LE_v = \varnothing$, but if $u$ and $v$ disagree on some element of $E_u \cap E_v$, there clearly cannot be an element in $X$ that contains both $u$ and $v$). Let $u, v \in \cL(Y)$ with shapes $F_u$ and $F_v$ such that $F_u \cap KF_v = \varnothing$. Since $u, v \in \cL(Y)$, we clearly have that $u, v \in \cL(X)$. We now show that $F_u \cap LF_v = \varnothing$. Indeed, since $K = L \cap H$, $F_u \cap KF_v = \varnothing$, and $F_u, F_v \subset H$,
        \begin{align*}
            F_u \cap LF_v &= \bigl(F_u \cap (L \cap H)F_v\bigr) \cup \bigl(F_u \cap (L \setminus H)F_v\bigr) \\
            &\subset \bigl(F_u \cap KF_v\bigr) \cup \bigl(H \cap (L \setminus H)H\bigr) \\
            &= \bigcup_{l \in L \setminus H} H \cap lH
        \end{align*}
        Since $H$ is a subgroup of $G$, for any $l \in L \setminus H$ we have $H \cap lH = \varnothing$, since $lH$ is a proper left coset of $H$. As such, $F_u \cap LF_v = \varnothing$, and so by the strong irreducibility of $X$, there exists $x \in X$ such that $x|_{F_u} = u$ and $x|_{F_v} = v$. Let $y = x|_H \in Y$, and so $y|_{F_u} = u$ and $y|_{F_v} = v$. Therefore, $Y$ is strongly irreducible.
    \end{proof}
\end{lemma}

Lastly, factor maps which are defined by block maps on a base shifts can be extended to a factor map of the free extension of the base shift in a natural manner. This does not apply to arbitrary factor maps from a free extension shifts however, and only applies to factor maps whose block maps are defined on a subset of the group for the base shift. In the result, for a function $\phi:X \to X$, we denote by $(\phi)^{H \backslash G}$ the product function on $X^{H \backslash G}$. If $X$ is a topological space with $\phi$ continuous, $(\phi)^{H \backslash G}$ is continuous on $X^{H \backslash G}$ endowed with the product topology.

\begin{lemma} \label{FE:factor_decomp}
    Let $G$ be a group and $F \Subset H \le G$. For finite alphabets $\cA$ and $\cB$, let $X$ be a $H$-shift of $\cA^H$, and let $\beta:\cL_F(X) \to \cB$ be a block map. Then $\phi^H_\beta$ is a factor map, so let $Y = \phi^H_\beta(X)$ be a $H$-shift of $\cB^H$. Let $\kappa_\cA$ denote the construction function on $\cA^{H : G}$ restricted to $X^{H : G}$ (and co-domain restricted to its image), and similarly let $\kappa_\cB$ denote the construction function on $\cB^{H:G}$ restricted to $Y^{H : G}$ using the same $\zeta$ as for $\kappa_\cA$. Then $\phi^G_\beta$ can be written as
    \[
        \phi^G_\beta = \kappa_\cB \circ (\phi^H_\beta)^{H \backslash G} \circ \kappa_\cA^{-1}.
    \]
    \begin{proof}
        Let $x \in X^{\uparrow G}$. Then for any $g \in G$, let $h \in H$ and $d \in C_{\kappa_\cA} = C_{\kappa_\cB}$ such that $g = hd$, and
        \begin{align*}
            \bigl((\kappa_\cB \circ (\phi^H_\beta)^{H \backslash G} \circ \kappa_\cA^{-1})(x)\bigr)(g) &= \left(\bigvee_{c \in C_{\kappa_\cB}} \sigma^{c^{-1}} \bigl(\phi^H_\beta(\kappa_\cA^{-1}(x)_c)\bigr)\right)(g) \\
            &= \Bigl(\sigma^{d^{-1}}\bigl(\phi^H_\beta(\kappa_\cA^{-1}(x)_d)\bigr)\Bigr)(hd) \\
            &= \phi^H_\beta\bigl((\sigma^dx)|_H\bigr)(h) \\
            &= \beta\Bigl(\bigl(\sigma^h\bigl((\sigma^dx)|_H\bigr)\bigr)|_F\Bigr) \\
            &= \beta\bigl(\bigl((\sigma^{hd}x)|_{Hh^{-1}}\bigr)|_F\bigr) \\
            &= \beta\bigl((\sigma^gx)|_F\bigr) \\
            &= \bigl(\phi^G_\beta(x)\bigr)(g)
        \end{align*}
    \end{proof}
\end{lemma}

A direct consequence of the previous result is that certain factors of a free extension are equal to the free extension of a factor of the base shift. This is not the case for all factors, however the property is quite useful nevertheless.

\begin{cor} \label{FE:factor_commutes}
    Let $G$ be a group, $H \le G$, $\cA$ and $\cB$ be finite alphabets, and $X$ be an $H$-shift of $\cA^H$. Let $F \Subset H$ and $\beta:\cL_F(X) \to \cB$ be a block map. Then
    \[
        \phi^G_\beta(X^{\uparrow G}) = \phi^H_\beta(X)^{\uparrow G}
    \]
    \begin{proof}
        By Lemma \ref{FE:factor_decomp}, $\phi^G_\beta = \kappa \circ (\phi^H_\beta)^{H \backslash G} \circ \kappa^{-1}$, and therefore,
        \begin{align*}
            \phi^G_\beta(X^{\uparrow G}) &= (\kappa \circ (\phi^H_\beta)^{H \backslash G} \circ \kappa^{-1})(X^{\uparrow G}) \\
            &= \bigl(\kappa \circ (\phi^H_\beta)^{H \backslash G}\bigr)\bigl(\kappa^{-1}(X^{\uparrow G})\bigr) \\
            &= \kappa \Bigl( (\phi^H_\beta)^{H \backslash G}\bigl(X^{H \backslash G}\bigr) \Bigr) \\
            &= \kappa \bigl( \phi^H_\beta(X)^{H \backslash G}\bigr) \\
            &= \phi^H_\beta(X)^{\uparrow G}.
        \end{align*}
    \end{proof}
\end{cor}

\subsection{Applications of free extensions to shifts on groups}

Using free extensions, it is possible to analyze shifts on arbitrary groups, though only to an extent. First, we can look at SFTs on arbitrary groups. We use the following result extensively in the study of SFTs on locally finite groups in particular, however it applies in full generality to all groups.

\begin{lemma} \label{SFT_FE}
    Let $G$ be a group, and $X$ a $G$-SFT. Then there exists $F \Subset G$ and $\langle F \rangle$-SFT $Y$ such that $X = Y^{\uparrow G}$. In other words, every SFT on a group $G$ is the free extension of an SFT on a finitely generated subgroup of $G$.
    \begin{proof}
        Since $X$ is an SFT, let $F \Subset G$ be a forbidden shape for $X$, so $X = \cX^G[\cF_F(X)]$. Let $H = \langle F \rangle \le G$, which makes $H$ finitely generated. Additionally, since $F \Subset H$, the $H$-shift $Y = \cX^H[\cF_F(X)]$ is an $H$-SFT. By Lemma \ref{FE:forbidden_ext},
        \[
            Y^{\uparrow G} = \cX^H[\cF_F(X)]^{\uparrow G} = \cX^G[\cF_F(X)] = X,
        \]
        which proves the desired result.
    \end{proof}
\end{lemma}

In the case that $G$ itself is finitely generated, it may be that $F \Subset G$ is such that $\langle F \rangle = G$, and so $X = Y = Y^{\uparrow G}$, which is a trivial result. In the case that $G$ is infinitely generated however, this can never be the case, and leads to interesting results such as the following.

\begin{cor} \label{inf_gen_entropies}
    Let $G$ be an infinitely generated amenable group. Then
    \[
        \cE(G) = \bigcup_{F \Subset G} \cE(\langle F \rangle).
    \]
    \begin{proof}
        This follows immediately from the previous lemma and Proposition \ref{FE:entropy}.
    \end{proof}
\end{cor}

In addition to SFTs, which are defined by a finite forbidden shape, strongly irreducible shifts are largely defined by the finite shape $K$, and using a similar technique as the lemma above, we can show that any strongly irreducible shift on a group is the free extension of a strongly irreducible shift on a finitely generated subgroup.

\begin{lemma} \label{SI_FE}
    Let $G$ be a group and $X$ be a strongly irreducible $G$-shift. Then there exists $F \Subset G$ and strongly irreducible $\langle F \rangle$-shift $Y$ such that $X = Y^{\uparrow G}$.
    \begin{proof}
        Since $X$ is strongly irreducible, let $K \Subset G$ be such that for any $u,v \in \cL(X)$ with shapes $F_u$ and $F_v$ respectively such that $F_u \cap KF_v = \varnothing$, then there exists $x \in X$ such that $x|_{F_u} = u$ and $x|_{F_v} = v$. We will now show that $K$ is an option for the finite set $F$ in the statement of the lemma. For any $F \Subset H$ let $\bF_F = \cF_F(X)$, and define $Y_F = \cX^H[\bF_F]$. Then since $F$ is finite, $Y_F$ is an $H$-SFT for each $F$. Furthermore, by Lemma \ref{FE:forbidden_ext}, we have that $Y_F^{\uparrow G} = (\cX^H[\bF_F])^{\uparrow G} = \cX^G[\bF_F]$, and so clearly $X \subset Y_F^{\uparrow G}$. As such, $X \subset \bigcap_{F \Subset H} Y_F^{\uparrow G}$, and by Lemma \ref{FE:intersect}, we have that
        \[
            \bigcap_{F \Subset H} Y_F^{\uparrow G} = \left(\bigcap_{F \Subset H} Y_F\right)^{\uparrow G}.
        \]
        Let $Y = \bigcap_{F \Subset H} Y_F$, which is an $H$-shift, so we have $X \subset Y^{\uparrow G}$. We now show that $Y^{\uparrow G} \subset X$.

        Let $\kappa$ be a construction function on $\cA^{H : G}$, let $z \in Y^{\uparrow G}$, and let $g \in G$ and $F \Subset H$. By the construction of $Y^{\uparrow G}$, we have $z \in Y_F^{\uparrow G} = \cX^G[\bF_F]$, and so $z|_F \notin \bF_F = \cF_F(X)$. Therefore, there exists $x_F \in X$ such that $z|_F = x_F|_F$. In particular, this shows that the set $E_F = [z|_F] \cap X$ is nonempty and closed. Additionally, since for each $g \in G$ we have $\sigma^gz \in Y^{\uparrow G}$, we also have that $[(\sigma^gz)|_F] \cap X = \sigma^g([z|_{Fg}] \cap X)$ is nonempty and closed. Since $\sigma^g$ is a homeomorphism on $\cA^G$ and $X$, we have that $E_{Fg} = [z|_{Fg}] \cap X$ is a nonempty closed subset of $X$. As such,
        \[
            \sG = \{E_{Fg} : F \Subset H, g \in G\}
        \]
        is a collection of nonempty closed subsets of $X$. We now show that $\sG$ has the finite intersection property.

        Let $E_{F_1g_1}, E_{F_2g_2}, \dots, E_{F_ng_n} \in \sG$. Note that since $F_i \Subset H$ and $g_i \in G$, we have that $F_ig_i \Subset Hc_i$ for some unique $c_i \in C_\kappa$. If we have that $F_jg_j \Subset Hc_i$ for some $j \ne i$, then $F_ig_i \cup F_jg_j \Subset Hc_i$, giving $(F_ig_i \cup F_jg_j)c_i^{-1} \Subset H$, and so we have
        \[
            E_{F_ig_i} \cap E_{F_jg_j} = [z|_{F_ig_i}] \cap [z|_{F_jg_j}] \cap X = [z|_{(F_ig_i \cup F_jg_j)c_i^{-1}c_i}] \cap X = E_{F_ig_i \cup F_jg_j},
        \]
        which is nonempty, and so we may assume without loss of generality that $c_i \ne c_j$ for $i \ne j$. For finite induction, we have that $E_{F_1g_1}$ is nonempty, so suppose that we have shown $\bigcap_{i=1}^k E_{F_ig_i}$ is nonempty for some $k < n$. As such, there exists an element $x \in X$ such that $u = x|_{\bigcup_{i=1}^k F_ig_i} = z|_{\bigcup_{i=1}^k F_ig_i}$, meaning that $u \in \cL(X)$. Let $v = z|_{F_{k+1}g_{k+1}}$, and note that $v \in \cL(X)$, as $E_{F_{k+1}g_{k+1}}$ is nonempty. Now, since $F_ig_i \subset Hc_i$ for each $i$, and $K \subset H$ by definition of $H$, we have that
        \[
            \left(\bigcup_{i=1}^k F_ig_i\right) \cap K(F_{k+1}g_{k+1}) \subset \left(\bigcup_{i=1}^k Hc_i\right) \cap H(Hc_{k+1}) = \left(\bigcup_{i=1}^k Hc_i\right) \cap Hc_{k+1}.
        \]
        Since $c_i \ne c_j$ for $i \ne j$, we have in the rightmost set an intersection of a right coset with a union of distinct right cosets, which is necessarily empty, and so we have
        \[
            \left(\bigcup_{i=1}^k F_ig_i\right) \cap K(F_{k+1}g_{k+1}) = \varnothing.
        \]
        By the strong irreducibility of $X$, there exists $x \in X$ such that $x|_{\bigcup_{i=1}^k F_ig_i} = u$ and $x|_{F_{k+1}g_{k+1}} = v$. This gives that $x \in \bigcap_{i=1}^{k+1} E_{F_ig_i}$, so this set is nonempty. By inducing until $n$, we obtain that $\bigcap_{i=1}^n E_{F_ig_i}$ is nonempty. As such, $\sG$ has the finite intersection property.

        Since $X$ is a closed subset of $\cA^G$, which is compact, we have that $X$ is compact. As such, since $\sG$ is a collection of closed subsets of $X$ with the finite intersection property, $\bigcap \sG$ is also nonempty, in particular there exists $x \in X$ such that
        \[
            x \in \bigcap_{g \in G} \bigcap_{F \Subset H} E_{Fg}.
        \]
        With $\{e\} \Subset H$, this gives that for each $g \in G$, $x \in E_{\{g\}} = [z|_{\{g\}}] \cap X$, which gives that $x(g) = z(g)$, and therefore $x = z$, which gives that $z \in X$. Since $z \in Y^{\uparrow G}$ was arbitrary, we have shown that $Y^{\uparrow G} \subset X$, and therefore $X = Y^{\uparrow G}$.

        Finally, by Lemma \ref{FE:si}, since $X = Y^{\uparrow G}$ is strongly irreducible, we have that $Y$ is strongly irreducible.
    \end{proof}
\end{lemma}

Further study into the properties of free extensions and which properties translate from a free extension to its base shift and vice-versa, may prove to show that the study of SFTs on arbitrary groups may be reducible to studying SFTs on finitely generated groups. While we do not require any further properties for the results of this paper, it may be fruitful to explore other such properties in the context of free extensions.

\section{Locally Finite Groups} \label{LF}

With the theory of free extensions sufficiently developed, we may move on to proving properties of SFTs on locally finite groups. This section contains all parts of the proofs of Theorems \ref{thmI} and \ref{thmII}.

We first begin by introducing the following construction, which applies to any group $G$ which is not locally finite, and which will be referenced throughout the remainder of the section.

\begin{defn}
    Let $G$ be a non-locally finite group, and $\cA = \{0, 1\}$. Since $G$ is non-locally finite, there exists an infinite, finitely generated group $H \le G$. Let $S \Subset H$ be such that $e \in S$ and $\langle S \rangle = H$. Then taking $\bF = \cA^S \setminus \{0^S, 1^S\}$, where $0^S$ and $1^S$ are the constant $0$ and $1$ patterns, let $\2_H = \cX^H[\bF]$.
    
    $\2_H$ is clearly an SFT from this construction, and in particular contains exactly 2 points, the constant 0 and 1 patterns on $H$, which will be denoted $0^H$ and $1^H$ respectively. By Lemma \ref{FE:sft}, $\2_H^{\uparrow G}$ is also an SFT.
\end{defn}

\subsection{Proof of Theorem \ref{thmI}}

We now have everything needed to prove Theorem \ref{thmI}. Each of the results in this section which contributes to the Theorem will be marked with the implication that it provides. For instance, the following result is marked as (\ref{thmI}\ref{I:lf}$\implies$\ref{thmI}\ref{I:fe}) to indicate that it provides the implication that if $G$ is locally finite, then every $G$-SFT is the free extension of an SFT on a finite subgroup of $G$. Many of these results follow readily from the properties of free extensions developed in the previous section. Theorem \ref{thmI} is restated below for convenience.

\begin{thm*}[I]
    Let $G$ be a group. Then the following are equivalent.fallow
    \begin{alphenum}
        \item $G$ is locally finite.
        \item Every $G$-SFT is the free extension of some SFT on a finite subgroup of $G$.
        \item Every $G$-SFT is strongly irreducible.
        \item Every strongly irreducible $G$-shift is a $G$-SFT.
        \item Every sofic $G$-shift is a $G$-SFT.
        \item For every $G$-SFT $X$, $\Aut(X)$ is locally finite.
    \end{alphenum}
\end{thm*}

We begin by proving the following chain of implications:
\begin{center}
    \ref{thmI}\ref{I:lf}$\implies$\ref{thmI}\ref{I:fe}$\implies$\ref{thmI}\ref{I:ir}$\implies$\ref{thmI}\ref{I:lf}
\end{center}
The first of these implications follows directly from Lemma \ref{SFT_FE}.

\begin{prop}[\ref{thmI}\ref{I:lf}$\implies$\ref{thmI}\ref{I:fe}] \label{LFG:SFT_are_ext}
    Let $G$ be a locally finite group, and $X$ a $G$-SFT. Then there exists $H \ll G$ and an $H$-SFT $Y$ such that $X = Y^{\uparrow G}$.
    \begin{proof}
        By Lemma \ref{SFT_FE}, there exists $F \Subset G$ and $\langle F \rangle$-SFT $Y$ such that $X = Y^{\uparrow G}$. But $G$ is locally finite, so $H = \langle F \rangle$ is finite, which gives the desired result.
    \end{proof}
\end{prop}

Next we show that if $G$ is a group for which every $G$-SFT is the free extension of a shift on a finite subgroup of $H$, then every $G$-SFT is strongly irreducible. In fact, we can show the following result,  which is stronger; if $X$ is a $G$-SFT for which there exists $H \ll G$ and $H$-SFT $Y$ such that $X = Y^{\uparrow G}$, then $X$ is strongly irreducible. 
\begin{lemma}[\ref{thmI}\ref{I:fe}$\implies$\ref{thmI}\ref{I:ir}]
    Let $G$ be a group and $X$ be a $G$-SFT such that there exists $H \ll G$ and $H$-SFT $Y$ such that $X = Y^{\uparrow G}$. Then $X$ is strongly irreducible.
    \begin{proof}
        Since $H$ is finite, $Y$ is vacuously strongly irreducible with $K = H$. By Lemma \ref{FE:si}, $X = Y^{\uparrow G}$ is strongly irreducible.
    \end{proof}
\end{lemma}

Lastly, we prove the final implication by contrapositive, where we use the SFT $\2^{\uparrow G}_H$ as an example of an SFT on non-locally finite groups which is not strongly irreducible.
\begin{lemma}[\ref{thmI}\ref{I:ir}$\implies$\ref{thmI}\ref{I:lf}]
    Let $G$ be a non-locally finite group. Then there exists a $G$-SFT which is not strongly irreducible.
    \begin{proof}
        Let $H \le G$ be an infinite, finitely generated subgroup of $G$, which must exist because $G$ is not locally finite.
    
        To show that $\2_H^{\uparrow G}$ is not strongly irreducible, it is necessary to show that for all $K \Subset G$, there exist patterns $u, v \in \cL(\2_H^{\uparrow G})$ with shapes $F_u$ and $F_v$ respectively such that $F_u \cap KF_v = \varnothing$, but there is no $x \in X$ with $x|_{F_u} = u$ and $x|_{F_v} = v$.
        
        Let $K \Subset G$. Since $K$ is finite, it must be that $H \setminus K$ is nonempty, so let $h \in H \setminus K$. Let $u = 0^{\{h\}}$ and $v = 1^{\{e\}}$, which are trivially in $\cL(\2_H^{\uparrow G})$. Then $F_u = \{h\}$ and $F_v = \{e\}$, and clearly since $h \notin K$, we have $F_u \cap KF_v = \{h\} \cap K = \varnothing$. But, for any $x \in X$,  it must be that $x|_H \in \{0^H, 1^H\}$, and therefore $x|_{\{h\}} = x|_{\{e\}}$, so it cannot be that $x|_{F_u} = u$ and $x|_{F_v} = v$ simultaneously.
        
        Therefore, $\2_H^{\uparrow G}$ is not strongly irreducible.
    \end{proof}
\end{lemma}

Next, we shall prove that \ref{thmI}\ref{I:lf}$\implies$\ref{thmI}\ref{I:ir2}, and prove the converse direction in the subsection immediately following, as we will need an example introduced then.

\begin{lemma}[\ref{thmI}\ref{I:lf}$\implies$\ref{thmI}\ref{I:ir2}]
    Let $G$ be a locally finite group, and $X$ a strongly irreducible $G$-shift. Then $X$ is a $G$-SFT.
    \begin{proof}
        By Lemma \ref{SI_FE}, there exists $F \Subset G$ and strongly irreducible $\langle F \rangle$-shift $Y$ such that $X = Y^{\uparrow G}$. Since $G$ is locally finite and $F$ is finite, $H = \langle F \rangle$ is finite, and therefore $Y$ is an $H$-SFT. By Lemma \ref{FE:sft}, $Y^{\uparrow G} = X$ is a $G$-SFT.
    \end{proof}
\end{lemma}

\subsubsection{Sofic shifts on locally finite groups}

Next, we prove the following implication involving the statement that every sofic $G$-shift is a $G$-SFT.
\begin{center}
    \ref{thmI}\ref{I:lf}$\iff$\ref{thmI}\ref{I:sof}
\end{center}

First, we show directly that all sofic $G$-shifts on locally finite groups are $G$-SFTs. 
\begin{lemma}[\ref{thmI}\ref{I:lf}$\implies$\ref{thmI}\ref{I:sof}]
    Let $G$ be a locally finite group, $X$ be an SFT, and $\phi$ be a factor map. Then $\phi(X)$ is an SFT.
    \begin{proof}
        By Proposition \ref{LFG:SFT_are_ext}, there exists $H \ll G$ and $H$-SFT $Y$ such that $X = Y^{\uparrow G}$. Let $F \Subset G$ and let $\beta:\cL_F(X) \to \cB$ be a block map such that $\phi = \phi^G_\beta$. Let $K = \langle H \cup F \rangle$, which is finite because $G$ is locally finite, and let $Z = Y^{\uparrow K}$. By Lemma \ref{FE:mult_ext}, we have $X = Y^{\uparrow G} = (Y^{\uparrow K})^{\uparrow G} = Z^{\uparrow G}$. Since $F \subset K$, Lemma \ref{FE:factor_commutes} gives that
        \[
            \phi(X) = \phi^G_\beta(Z^{\uparrow G}) = \phi^K_\beta(Z)^{\uparrow G}.
        \]
        Since $K$ is finite, the $K$-shift $\phi^K_\beta(Z) \subset \cA^K$ is an SFT, and by Lemma \ref{FE:sft}, we obtain that $\phi^K_\beta(Z)^{\uparrow G} = \phi(X)$ is an SFT.
    \end{proof}
\end{lemma}

For the converse result, we will prove the contrapositive by constructing for any non-locally finite group, a sofic shift which is not an SFT. We begin with the construction.

\begin{defn}[Example 1.11 \cite{barbieri_thesis}]
    Let $H$ be an infinite, finitely generated group, and let $S \Subset H$ such that $S = S^{-1}$, $e \notin S$, and $H = \langle S \rangle$. The \textit{even $H$-shift} $S_{\textrm{even}} \subset \{0,1\}^H$ is the set of all full patterns $x$ such that any maximal finite connected component of $x^{-1}(1) \subset H$ in the Cayley graph $\Gamma(H, S)$ has even size. Said otherwise, each finite connected component of ones has even size.
\end{defn}

Proposition 1.10 of \cite{barbieri_thesis} gives that $S_{\textrm{even}}$ is a sofic $H$-shift, but not an $H$-SFT. Using this, we can prove the converse result.

\begin{lemma}[\ref{thmI}\ref{I:sof}$\implies$\ref{thmI}\ref{I:lf}]
    Let $G$ be a non-locally finite group. Then there exists a sofic $G$-shift which is not a $G$-SFT.
    \begin{proof}
        Let $H \le G$ be infinite and finitely generated. Then $S_{\textrm{even}}$ as defined above is a sofic $H$-shift, but not an $H$-SFT. Let $X$ be an $H$-SFT and $\phi:X \to S_{\textrm{even}}$ be a factor map, which must exist by the soficity of $S_{\textrm{even}}$. Then there exists $F \Subset H$ and a block map $\beta:\cL_F(X) \to \{0, 1\}$ such that $\phi = \phi^H_\beta$. Then by Lemma \ref{FE:sft}, it follows $X^{\uparrow G}$ is a $G$-SFT, and by Lemma \ref{FE:factor_commutes}, we obtain
        \[
            S_{\mathrm{even}}^{\uparrow G} = \phi^H_\beta(X)^{\uparrow G} = \phi^G_\beta(X^{\uparrow G}),
        \]
        and therefore $S_{\mathrm{even}}^{\uparrow G}$ is sofic. By the contrapositive of Lemma \ref{FE:sft} however, $S_{\textrm{even}}^{\uparrow G}$ is not an SFT.
    \end{proof}
\end{lemma}

In addition to $S_{\textrm{even}}$ being a sofic $H$-shift, we also have that it is strongly irreducible. With $K = (S \cup \{e\})^2$, and two patterns $u, v \in S_{\textrm{even}}$ with shapes $F_u$ and $F_v$ such that $F_u \cap KF_v = \varnothing$, we may extend $u$ to a pattern on $(S \cup \{e\})F_u$ by using the symbol 0 or 1 in a manner that ensures this extension has an even number of 1s in any connected component of 1s. The same can be done for $v$. By the definition of $K$, we have that these extensions are disjoint, so let $x \in \{0,1\}^H$ which matches these extensions of $u$ and $v$, and is 0 elsewhere. Since the extensions of $u$ and $v$ have connected components of 1s of even size, $x$ only has connected components of 1s of even size, even if a connected component in the extension of $u$ is connected with a connected component of the extension of $v$, since both individually have even size. As such, we have the following result.

\begin{lemma}[\ref{thmI}\ref{I:ir2}$\implies$\ref{thmI}\ref{I:lf}]
    Let $G$ be a non-locally finite group. Then there exists a strongly irreducible $G$-shift which is not a $G$-SFT.
    \begin{proof}
        Let $H \le G$ be infinite and finitely generated. Then $S_{\textrm{even}}$ as defined above is a strongly irreducible $H$-shift, but not an $H$-SFT. By Lemma \ref{FE:sft}, $S_{\textrm{even}}^{\uparrow G}$ is not a $G$-SFT, and by Lemma \ref{FE:si}, $S_{\textrm{even}}^{\uparrow G}$ is strongly irreducible.
    \end{proof}
\end{lemma}

\subsubsection{Automorphism groups for locally finite SFTs}

Finally, we prove the last implications for Theorem \ref{thmI} in the following manner.
\begin{center}
    \ref{thmI}\ref{I:lf}$\iff$\ref{thmI}\ref{I:aut} 
\end{center}

First, we show that the automorphism group for any SFT on a locally finite groups is locally finite.
\begin{lemma}[\ref{thmI}\ref{I:lf}$\implies$\ref{thmI}\ref{I:aut}]
    Let $G$ be locally finite and $X$ a $G$-SFT. Then $\Aut(X)$ is locally finite.
    \begin{proof}
        Let $F \Subset G$ be a forbidden shape for $X$ so that $X = \cX^G[\cF_F(X)]$. Let $E = \{\phi_1, \phi_2, \dots, \phi_n\} \subset \Aut(X)$ be a finite set of autoconjugacies, and let $K = \langle E \rangle$. Without loss of generality, $E$ may be assumed to be symmetric. Then for each $\phi_i$, there exists $F_i \Subset G$ and block maps $\beta_i:\cL_{F_i}(X) \to \cA$ such that $\phi_i = \phi^G_{\beta_i}$. Now, let
        \[
            H = \left\langle F \cup \bigcup_{i = 1}^n F_i \right\rangle.
        \]
        $H$ must be finite, since $G$ is locally finite. Then, since $F \subset H$, it is the case that $X = \cX^G[\cF_H(X)]$, and by Lemma \ref{FE:forbidden_ext}, we have $X = \cX^G[\cF_H(X)] = \cX^H[\cF_H(X)]^{\uparrow G}$. For simplicity, let $Y = \cX^H[\cF_H(X)]$. Additionally, by Corollary \ref{FE:factor_commutes}, for each $i$, we obtain
        \[
            Y^{\uparrow G} = \phi_i(Y^{\uparrow G}) = \phi^H_{\beta_i}(Y)^{\uparrow G},
        \]
        and therefore $Y = \phi^H_{\beta_i}(Y)$, which gives $\phi^H_{\beta_i} \in \Aut(Y)$. As such, let
        \[
            \hat{E} = \{\phi^H_{\beta_i} : i \in [1, n]\},
        \]
        and let $\hat{K} = \langle \hat{E} \rangle \le \Aut(Y)$.
        
        Now, let $\kappa$ be a construction function on $\cA^{H:G}$, and define a map $\gamma:\hat{K} \to \Aut(X)$ by
        \[
            \gamma(\psi) = \kappa \circ (\psi)^{H \backslash G} \circ \kappa^{-1},
        \]
        where $\kappa$ is taken to be restricted to $\kappa:Y^{H : G} \to X$. Then for any $\phi \in K$, there exists $k$ and $i_1, i_2, \dots, i_k$ such that $\phi = \phi_{i_1} \circ \phi_{i_2} \circ \cdots \circ \phi_{i_k}$. Let
        \[
            \psi = \phi^H_{\beta_{i_1}} \circ \phi^H_{\beta_{i_2}} \circ \cdots \circ \phi^H_{\beta_{i_k}}.
        \]
        Then clearly $\psi \in \hat{K}$, and by Lemma \ref{FE:factor_decomp}, we get
        \begin{align*}
            \gamma(\psi) &= \kappa \circ (\psi)^{H \backslash G} \circ \kappa^{-1} \\
            &= \kappa \circ \left(\phi^H_{\beta_{i_1}} \circ \phi^H_{\beta_{i_2}} \circ \cdots \circ \phi^H_{\beta_{i_k}}\right)^{H \backslash G} \circ \kappa^{-1} \\
            &= \kappa \circ \left(\phi^H_{\beta_{i_1}}\right)^{H \backslash G} \circ \kappa^{-1} \circ \kappa \circ \left(\phi^H_{\beta_{i_2}} \circ \cdots \circ \phi^H_{\beta_{i_k}}\right)^{H \backslash G} \circ \kappa^{-1} \\
            &= \phi^G_{\beta_{i_1}} \circ \kappa \circ \left(\phi^H_{\beta_{i_2}} \circ \cdots \circ \phi^H_{\beta_{i_k}}\right)^{H \backslash G} \circ \kappa^{-1} \\
            &= \phi_{i_1} \circ \phi^G_{\beta_{i_2}} \circ \kappa \circ \left(\phi^H_{\beta_{i_3}} \circ \cdots \circ \phi^H_{\beta_{i_k}}\right)^{H \backslash G} \circ \kappa^{-1} \\
            &= \phi_{i_1} \circ \phi_{i_2} \circ \cdots \circ \phi_{i_k} \\
            &= \phi,
        \end{align*}
        and therefore $K \subset \gamma(\hat{K})$, so $|K| \le |\gamma(\hat{K})| \le |\hat{K}| \le |\Aut(Y)|$. Since $H$ is finite and $Y$ is an $H$-shift, it must be that $\Aut(Y)$ is finite, and therefore $K$ must also be finite. Since $E$ was arbitrary, this gives that $\Aut(X)$ is locally finite.
    \end{proof}
\end{lemma}

Lastly, we show that if the automorphism group of an SFT is locally finite, then the group on which the SFT is defined must be locally finite.

\begin{lemma}[\ref{thmI}\ref{I:aut}$\implies$\ref{thmI}\ref{I:lf}]
    Let $G$ be a group. If for every $G$-SFT $X$, $\Aut(X)$ is locally finite, then $G$ is locally finite.
    \begin{proof}
        If every $G$-SFT $X$ satisfies $\Aut(X)$ is locally finite, in particular this is true of the full $G$-shift $\Sigma$ (on at least 2 symbols). Clearly, the map $\psi:G \to \Aut(\Sigma)$ defined by $\psi(g) = \sigma^g$ is an injective homomorphism, since for any $h \ne g$, we have $\sigma^h \ne \sigma^g$ on $\Sigma$, since it is possible to describe a configuration which gets sent to different configurations under $\sigma^h$ and $\sigma^g$. As such, $\psi(G) \le \Aut(\Sigma)$. Since $\Aut(\Sigma)$ is locally finite, $\psi(G)$ is locally finite. But $\psi(G)$ and $G$ are isomorphic, and therefore $G$ is locally finite.
    \end{proof}
\end{lemma}

\subsection{Proof of Theorem \ref{thmII}}

Next, we will prove Theorem \ref{thmII}. As with the previous section, results pertaining to certain implications in Theorem \ref{thmII} are marked. The main additional assumption we will need is that $G$ is a countable amenable group, rather than any group. Most of these results also depend heavily on the properties of free extensions developed in the previous section. We restate Theorem \ref{thmII} below for convenience.

\begin{thm*}[II]
    Let $G$ be a countable amenable group. Then the following are equivalent.
    \begin{alphenum}
        \item $G$ is locally finite 
        \item If $X$ is a nonempty $G$-SFT with $h(X) = 0$, then $X = \{x\}$, where $x$ is a fixed point.
        \item Every $G$-SFT is entropy minimal.
        \item $G$ is locally typical and
        \[
            \cE(G) = \left\{\frac{\log(n)}{|H|} : H \ll G, n \in \N \right\} \subset \Q^+_{\log}.
        \]
        \item Every $G$-SFT has a unique measure of maximal entropy.
    \end{alphenum}
\end{thm*}

Each of the equivalences in the theorem will be shown individually to be equivalent to \ref{thmII}\ref{II:lf}. We begin by showing its equivalence to \ref{thmII}\ref{II:0_ent}. Additionally, note that all countable locally finite groups are amenable, so we omit amenable as an assumption for a few of the results.

\subsubsection{Zero entropy SFTs on locally finite groups}

We begin by showing that zero entropy SFTs on locally finite groups consist of single fixed points.
\begin{lemma}[\ref{thmII}\ref{II:lf}$\implies$\ref{thmII}\ref{II:0_ent}]
    Let $G$ be a countable locally finite group. Then if $X$ is a non-empty $G$-SFT with $h(X) = 0$, then $X = \{x\}$ for some fixed point $x$.
    \begin{proof}
        Let $X$ be a $G$-SFT with $h(X) = 0$. Then by assumption, $X = Y^{\uparrow G}$ for some $H \ll G$ and $H$-shift $Y$. By Proposition \ref{FE:entropy}, we have $h(X) = h(Y) = 0$. Since $H$ is finite,
        \[
            0 = h(Y) = \frac{1}{|H|}\log(|Y|), 
        \]
        which implies that $|Y| = 1$. Then $|Y^{H \backslash G}| = 1$, and therefore $|X| = |\kappa(Y^{H \backslash G})| = 1$ for any construction function $\kappa$, so $X = \{x\}$ for the only $x \in X$. Since $X$ is shift invariant, it must be that $x$ is a fixed point.
    \end{proof}
\end{lemma}

To show the converse, recall the definition of the SFT $\2_H$ from the beginning of Section \ref{LF}.
\begin{lemma}[\ref{thmII}\ref{II:0_ent}$\implies$\ref{thmII}\ref{II:lf}] \label{2_0_ent}
    Let $G$ be a countable amenable non-locally finite group. Then there exists a $G$-SFT $X$ with 0 topological entropy, however $|X| > 1$.
    \begin{proof}
        Let $H \le G$ be an infinite, finitely generated subgroup. By Proposition \ref{FE:entropy}, we have $h(\2_H^{\uparrow G}) = h(\2_H)$. Since $G$ is countable and amenable, and $H \le G$, it is the case that $H$ is also countable and amenable, so let $\{F_n\}_{n=1}^\infty$ be a F{\o}lner sequence for $H$. Since $\2_H$ contains exactly 2 points, $0^H$ and $1^H$, it is clear to see that $\cL_{F_n}(\2_H) = \{0^{F_n}, 1^{F_n}\}$, and therefore $|\cL_{F_n}(\2_H)| = 2$. Additionally, it must be that $\lim_{n \to \infty} |F_n| = \infty$, because $H$ is an infinite subgroup. Then
        \[
            h(\2_H^{\uparrow G}) = h(\2_H) = \lim_{n \to \infty} |F_n|^{-1} \log\bigl(|\cL_{F_n}(\2_H)|\bigr) = \log(2) \lim_{n \to \infty} |F_n|^{-1} = 0.
        \]
        Also, since $|\2_H| = 2$, $|\2_H^{\uparrow G}| > 1$, which gives the desired result.
    \end{proof}
\end{lemma}

\subsubsection{Entropy minimality of SFTs on locally finite groups}

Recall that a $G$-SFT $X$ is entropy minimal if for every $G$-shift $Y \subsetneq X$, we have $h(Y) < h(X)$. The following result shows that for a countable locally finite group, every SFT on the group is entropy minimal.
\begin{lemma}[\ref{thmII}\ref{II:lf}$\implies$\ref{thmII}\ref{II:ent_min}]
    Let $G$ be a countable locally finite group, and $X$ be a $G$-SFT. Then $X$ is entropy minimal.
    \begin{proof}
        Since $X$ is an SFT, let $F \Subset G$ be such that $X = \cX^G[\cF_F(X)]$. Let $Y \subsetneq X$ also be an SFT, and let $E \Subset G$ be such that $Y = \cX^G[\cF_E(Y)]$. Let $H = \langle F \cup E\rangle$, which is finite because $E$ and $F$ are finite and $G$ is locally finite. Also we have $E, F \subset H$, and therefore $X = \cX^G[\cF_H(X)]$ and $Y = \cX^G[\cF_H(Y)]$. By Lemma \ref{FE:forbidden_ext}, we obtain $X = (\cX^H[\cF_H(X)])^{\uparrow G}$ and $Y = (\cX^H[\cF_H(Y)])^{\uparrow G}$. Given that $Y \subsetneq X$, it must then be that $\cX^H[\cF_H(Y)] \subsetneq \cX^H[\cF_H(X)]$. With Proposition \ref{FE:sft}, this gives
        \begin{align*}
            h(Y) &= h(\cX^H[\cF_H(Y)]) = |H|^{-1}\log(|\cX^H[\cF_H(Y)]|) \\
            &< |H|^{-1}\log(|\cX^H[\cF_H(X)]|) = h(\cX^H[\cF_H(X)]) = h(X),
        \end{align*}
        where the strict inequality follows from the fact that $\log$ is strictly increasing, and this implies that $X$ is SFT-entropy minimal. Since $X$ is an SFT, and SFT-entropy minimality and entropy minimality are equivalent for SFTs, we have that $X$ is entropy minimal.
    \end{proof}
\end{lemma}

For the converse result about entropy minimality, we again use the SFT $\2_H$.
\begin{lemma}[\ref{thmII}\ref{II:ent_min}$\implies$\ref{thmII}\ref{II:lf}]
    Let $G$ be a countable amenable non-locally finite group. Then there exists a $G$-SFT $X$ which is not entropy minimal.
    \begin{proof}
        Let $H \le G$ be an infinite, finitely generated subgroup. We have that $\2_H^{\uparrow G}$ is a $G$-SFT, and $h(\2_H^{\uparrow G}) = 0$, an argument for which can be found in Lemma \ref{2_0_ent}. It is clear that $\{0^G\} \subset \2_H^{\uparrow G}$, and $\{0^G\}$ is clearly a $G$-shift, as it is conjugate to the full $G$-shift on 1 symbol. Additionally, $h(\{0^G\}) = 0$, and therefore $\2_H^{\uparrow G}$ is not entropy minimal.
    \end{proof}
\end{lemma}

\subsubsection{The set of SFT entropies for locally finite groups} \label{LF-sft_ent}

Next, we establish the set of all entropies that SFTs can obtain for locally finite groups. The following result shows that \ref{thmII}\ref{II:lf} implies \ref{thmII}\ref{II:sft_ent}. The first part of the implication is trivial; if $G$ is locally finite, then every finitely generated subgroup is finite, and therefore $G$ is locally typical. The second part of the implication is given below.
\begin{lemma}[\ref{thmII}\ref{II:lf}$\implies$\ref{thmII}\ref{II:sft_ent}] \label{ffe_entropies}
    Let $G$ be a countable locally finite group. Then
    \[
        \cE(G) = \left\{\frac{\log(n)}{|H|} : H \ll G, n \in \N \right\} \subset \Q^+_{\log}.
    \]
    \begin{proof}
        First, consider the case when $G$ is finite. Let $X$ be a $G$-SFT. Then $h(X) = \frac{\log(|X|)}{|G|} \in \Q^+_{\log}$, and so
        \[
            \cE(G) \subset \left\{\frac{\log(n)}{|H|} : H \ll G, n \in \N\right\},
        \]
        since $G \ll G$. Now let $H \ll G$ and $n \in \N$. Since $G$ and $H$ are finite, let $m = \frac{|G|}{|H|} \in \N$. Let $\cA$ be a finite alphabet with $|\cA| = n^m$. Then, let $X = \{a^G : a \in \cA\}$, which is a $G$-SFT, and $|X| = n^m$. Then
        \[
            h(X) = \frac{\log(|X|)}{|G|} = \frac{\log(n^m)}{|G|} = \frac{m\log(n)}{|G|} = \frac{\log(n)}{|H|},
        \]
        and therefore,
        \[
            \left\{\frac{\log(n)}{|H|} : H \ll G, n \in \N \right\} \subset \cE(G),
        \]
        which gives the desired result.
        
        If $G$ is infinite, then $G$ must be infinitely generated, and so by Corollary \ref{inf_gen_entropies},
        \[
            \cE(G) = \bigcup_{F \Subset G} \cE(\langle F \rangle).
        \]
        Since $G$ is locally finite, $H \ll G$ if and only if $H$ is finitely generated, which gives
        \begin{align*}
            \cE(G) &= \bigcup_{H \ll G} \cE(H) = \bigcup_{H \ll G} \left\{\frac{\log(n)}{|K|} : K \ll H, n \in \N\right\} \\
            &= \left\{\frac{\log(n)}{|H|} : H \ll G, n \in \N\right\}.
        \end{align*}
    \end{proof}
\end{lemma}

Many locally finite groups do not satisfy $\cE(G) = \Q^+_{\log}$, due to the lack of subgroups of certain orders. For example, $\bigoplus_{n \in \N} \Z/2\Z$ is locally finite, but only has subgroups of order $2^n$. There are locally finite groups which do attain $\cE(G) = \Q^+_{\log}$ however, with the most prominent example likely being Hall's universal group $\U$ \cite{hall}, which has the property that every countable locally finite group can be embedded within it, which includes all finite groups. As such, it has finite subgroups of every order, and so $\cE(\U) = \Q^+_{\log}$.

A direct converse of the previous lemma has been elusive to the author, which is the reason for the additional statement that $G$ is locally typical in \ref{thmII}\ref{II:sft_ent}. The following lemma gives the most general form of a converse that has been found by the author.
\begin{lemma} \label{LF:ent_imp_per}
    Let $G$ be a countable amenable group such that $\cE(G) \subset \Q^+_{\log}$. Then $G$ is periodic.
   
    \begin{proof}
        We proceed by the contrapositive. Let $G$ be a group which is not periodic, meaning there exists $h \in G$ whose order is infinite. Let $H = \langle h \rangle$ so that $H$ is isomorphic to $\Z$, and define $\bF = \{1^{\{e, h\}}\} \subset \{0,1\}^{\{e,h\}}$, and let $X = \cX^H[\bF]$. Since $G$ is amenable, and $H \le G$, it must be that $H$ is amenable. Then $X$ is conjugate to the well known golden mean shift on $\Z$, so $h(X) = \log(\varphi)$, where $\varphi = \frac{1 + \sqrt{5}}{2}$ is the golden ratio. $X$ is also clearly an SFT, so by Lemma \ref{FE:sft} the $G$-shift $X^{\uparrow G}$ is an SFT, and by Proposition \ref{FE:entropy} we have $h(X^{\uparrow G}) = h(X) = \log(\varphi)$. It is an elementary number theory exercise to show that $\varphi^n$ is irrational for all $n \in \N$, and so it must be that for any $n, m \in \N$, we have $\varphi^m \ne n$. Therefore $\forall n, m \in \N$, it is the case that $\log(\varphi) \ne \frac{\log(n)}{m}$, so $\log(\varphi) \notin \Q^+_{\log}$. But $\log(\varphi) \in \cE(G)$, and therefore $\cE(G) \not\subset \Q^+_{\log}$.
    \end{proof}
\end{lemma}

It remains to show that periodic but not locally finite groups have SFTs with entropy outside of $\Q^+_{\log}$, however it is in general quite difficult to construct SFTs on such groups in a manner conducive to computing its topological entropy. As a result, we instead add the statement that $G$ is locally typical, which removes the need to consider such groups.
\begin{lemma}[\ref{thmII}\ref{II:sft_ent}$\implies$\ref{thmII}\ref{II:lf}]
    Let $G$ be a countable amenable group which is locally typical, and $\cE(G) \subset \Q^+_{\log}$. Then $G$ is locally finite.
    \begin{proof}
        By Lemma \ref{LF:ent_imp_per}, $G$ is periodic. Let $F \Subset G$, and consider $H = \langle F \rangle$. Since $G$ is periodic, $H$ is periodic. Since $G$ is locally typical, $H$ is finite or not periodic, and therefore $H$ must be finite. Since $F \Subset G$ was arbitrary, $G$ is locally finite.
    \end{proof}
\end{lemma}

The author suspects that if $\cE(G) \subset \Q^+_{\log}$, then $\cE(G)$ must be locally finite. This would allow for \ref{thmII}\ref{II:sft_ent} to have the locally typical assumption removed, and only leave $\cE(G) \subset \Q^+_{\log}$.

\subsubsection{Measures of maximal entropy for SFTs on locally finite groups}

Finally, we show that every SFT on a countable locally finite group has a unique measure of maximal entropy, and that if every SFT on a countable amenable group has a unique measure of maximal entropy, then the group must be locally finite. First, we require a simple but powerful result about the topological structure of SFTs on countable locally finite groups.

\begin{lemma} \label{LFG:sft_top_basis}
    Let $G$ be a countable locally finite group, and let $X$ be a $G$-SFT. Then there exists a sequence $\{H_n\}_{n=1}^\infty$ with $H_n \le H_{n+1} \ll G$ for all $n$, such that $G = \bigcup_{n \in \N} H_n$, and there exist $H_n$-SFTs $Y_n$ such that $X = Y_n^{\uparrow G}$ for all $n$. Furthermore, the set
    \[
        \fB[\{Y_n\}] = \{[y] \cap X : n \in \N, y \in Y_n\}
    \]
    is a basis for the subspace topology on $X$.
    \begin{proof}
        First, since $X$ is a $G$-SFT, by Proposition \ref{LFG:SFT_are_ext}, there exists $H_1 \ll G$ and $H_1$-SFT $Y_1$ such that $X = Y_1^{\uparrow G}$. Then, since $G$ is countable, let $G = \{g_n : n \in \N\}$ be an enumeration of $G$. Define for $n \ge 2$,
        \[
            H_n = \langle H_1 \cup \{g_i : i < n\} \rangle.
        \]
        Since $H_1$ and $\{g_i : i < n\}$ are both finite, $H_n$ is finitely generated, and therefore finite. Furthermore, for any $g \in G$, there is some $n \in \N$ for which $g = g_n$, and clearly $g_n \in H_{n+1}$. Also, $H_n \le H_{n+1}$.
        
        Now, for each $n \ge 2$, let $Y_n = Y_1^{\uparrow H_n}$. By Lemma \ref{FE:mult_ext}, we obtain $X = Y_1^{\uparrow G} = (Y_1^{\uparrow H_n})^{\uparrow G} = Y_n^{\uparrow G}$.
        
        Finally, let $\fB$ be the standard basis of all cylinder sets for $X$. To show that $\fB[\{Y_n\}]$ is a basis for the topology on $X$, first note that $\fB[\{Y_n\}] \subset \fB$, and therefore it suffices to show that any set in $\fB$ can be constructed by sets in $\fB[\{Y_n\}]$. Let $w \in \cL(X)$ so that $[w] \cap X$ is nonempty, and let $F$ be the shape of $w$. Since $G = \bigcup_{n \in \N} H_n$ and $H_n \le H_{n+1}$, it follows there exists $N \in \N$ such that $F \subset H_N$. Then, it is clear that
        \[
            [w] \cap X = \bigcap_{z \in [w]_{H_N} \cap Y_N} [z] \cap X,
        \]
        where this intersection is finite, which implies that $\tau(\fB)$, the topology generated by $\fB$, is contained in $\tau(\fB[\{Y_n\}])$, so $\fB[\{Y_n\}]$ is a basis for the topology on $X$.
    \end{proof}
\end{lemma}

\begin{lemma}[\ref{thmII}\ref{II:lf}$\implies$\ref{thmII}\ref{II:mme}]
    Let $G$ be a countable locally finite group. Then for any $G$-SFT $X$, there exists a unique measure of maximal entropy.
    \begin{proof}
        Measures of maximal entropy exist for all shift actions of countable amenable groups, so let $\mu \in \cM(X)$ such that $h_\mu(X) = h(X)$.
        
        By Lemma \ref{LFG:sft_top_basis}, there exists $\{H_n\}_{n=1}^\infty$ and $H_n$-SFTs $Y_n$ such that $\fB[\{Y_n\}]$ is a basis for the topology on $X$, and therefore also generates the Borel $\sigma$-algebra on $X$. Furthermore, since $X = Y_n^{\uparrow G}$, Lemma \ref{FE:entropy} gives that
        \[
            h(X) = h(Y_1) = h(Y_n) = \frac{\log(|Y_1|)}{|H_1|} = \frac{\log(|Y_n|)}{|H_n|}
        \]
        for all $n \in \N$. Also note that $\{H_n\}$ is a F{\o}lner sequence for $G$, and therefore
        \[
            h_\mu(X) = \inf_n \frac{H_\mu(X, H_n)}{|H_n|}.
        \]
        As such, we obtain $h_\mu(X) \le \frac{H_\mu(X, H_n)}{|H_n|}$ for all $n \in \N$. But
        \[
            H_\nu(X, H_n) \le \log(|\cL_{H_n}(X)|) = \log(|Y_n|)\]
        for any $\nu \in \cM(X)$, and therefore
        \[
            \frac{\log(|Y_n|)}{|H_n|} = h(X) = h_\mu(X) \le \frac{H_\mu(X, H_n)}{|H_n|} \le \frac{\log(|Y_n|)}{|H_n|},
        \]
        so all of these quantities must be equal, which further implies that for every $n$ and $y \in Y_n$, we have $\mu[y] = \frac{1}{|Y_n|}$. This is true for any $n \in \N$, and therefore any measure of maximal entropy must take these specific values for every element of $\fB[\{Y_n\}]$. By the Carath\'eodory Extension Theorem, there exists a unique Borel probability measure with these properties, and therefore there exists only one measure of maximal entropy.
    \end{proof}
\end{lemma}

Though the previous proof does not explicitly mention how to construct the measure of maximal entropy, its construction is fairly simple. For a countable locally finite group $G$ and $G$-SFT $X$, take some $H \ll G$ and $H$-SFT $Y$ such that $X = Y^{\uparrow G}$. Let $\nu$ be a measure on $Y$ defined by $\nu(y) = \frac{1}{|Y|}$ for all $y \in Y$. Then $\mu = (\nu)^{H \backslash G} \circ \kappa^{-1}$ is an invariant measure of maximal entropy for $X$. Informally, $\mu$ is the uniform measure on $X$, which is obtained as the push forward of a product measure under the construction function. It can also be shown that $\mu$ is independent of choice of $H$ and $Y$ for which $X = Y^{\uparrow G}$.

For the converse result, we give an SFT on any non-locally finite group which has multiple measures of maximal entropy.

\begin{lemma}[\ref{II:mme}$\implies$\ref{II:lf}]
    Let $G$ be a countable amenable non-locally finite group. Then there exists a $G$-SFT $X$ which has multiple measures of maximal entropy.
    \begin{proof}
        Since $h(\2_H^{\uparrow G}) = 0$, the Variational Principle gives that for all $\mu \in \cM(X)$, we have $0 \le h_\mu(\2_H^{\uparrow G}) \le h(\2_H^{\uparrow G}) = 0$, and so $h_\mu(\2_H^{\uparrow G}) = h(\2_H^{\uparrow G})$. This means every measure $\mu \in \cM(X)$ is a measure of maximal entropy.
        
        Since $0^G$ and $1^G$ are both elements of $\2_H^{\uparrow G}$, the two Dirac measures $\delta_{0^G}$ and $\delta_{1^G}$ are distinct, and since both $0^G$ and $1^G$ are fixed points they are both invariant, and therefore contained within $\cM(\2_H^{\uparrow G})$. As such, $\2_H^{\uparrow G}$ has at least 2 measures of maximal entropy.
    \end{proof}
\end{lemma}

\section{Final Remarks} \label{FR}

The main results of this paper gives that the class of locally finite groups presents interesting dynamical behaviors that are unexpected in general. This combined with the converse results, which show these interesting behaviors are unique to locally finite groups, gives insights into the types of groups where interesting behavior is possible. As mentioned in Section \ref{INT}, Theorem \ref{thmII}\ref{II:0_ent} gives that the only groups for which there are only trivial zero-entropy dynamics are precisely the locally finite groups, and this property indirectly answers in the affirmative Question 3.19 of Barbieri \cite{barbieri}: \dquote{Does there exist an amenable group $G$ and a $G$-SFT which does not contain a zero-entropy $G$-SFT?} For any countable locally finite group $G$, take any finite $H \ll G$ with $|H| > 1$, and pick any $H$-SFT $Y$ which does not contain any fixed points. Then $X = Y^{\uparrow G}$ also does not contain a fixed point, and therefore contains no zero-entropy SFTs. This answer to the question leads to the following refinement of the question, as infinite locally finite groups are necessarily infinitely generated:

\begin{ques}
    Does there exist an infinite, \textit{finitely generated} amenable group $G$ and a $G$-SFT which does not contain a zero-entropy $G$-SFT?
\end{ques}

Theorem \ref{thmII}\ref{II:sft_ent} also aids in the overall classification of the possible sets which are attainable as the set of entropies of SFTs on a specific group. In the case that $G$ is locally finite, $\cE(G) \subset \Q^+_{\log}$ (and in particular, an exact form for $\cE(G)$ is known). In the case that $G$ is not periodic, then it contains an element of infinite order (and therefore a subgroup isomorphic to $\Z$), and thus by Lemma \ref{SFT_FE}, $\cE(\Z) \subset \cE(G)$. Though more research is needed to classify $\cE(G)$ exactly for these types of groups (such as the work of Barbieri \cite{barbieri}), at least it is known that $\Z$-SFT entropies are attainable. The remaining class of groups are the infinite, periodic, but not locally finite groups. We have shown in Lemma \ref{LF:ent_imp_per} that $\cE(G) \subset \Q^+_{\log}$ does imply that the $G$ is periodic, however it is unclear whether the following can be answered in the affirmative:
\begin{ques}
    If $G$ is a countable amenable group such that $\cE(G) \subset \Q^+_{\log}$, then must it be the case that $G$ is locally finite? If not, is
    \[
        \cE(G) = \left\{\frac{\log(n)}{|H|} : H \ll G, n \in \N\right\}
    \]
    sufficient to conclude that $G$ is locally finite?
\end{ques}
Answering either of these questions in the affirmative would permit the locally typical statement in \ref{thmII}\ref{II:sft_ent} to be dropped, leading to a strictly stronger result. Following the method used in proving that $\cE(G) \subset \Q^+_{\log}$ implies periodicity, it would suffice to produce for any infinite, finitely generated, periodic group $H$, an $H$-SFT with entropy outside of $\Q^+_{\log}$. Then, for any infinite periodic group $G$ which is not locally finite, it must contain an infinite, finitely generated, periodic subgroup $H$ (potentially the whole group), and this SFT can be defined on $H$, and then freely extended to $G$ with the same entropy. Defining SFTs is not difficult in general; the primary difficulty is in computing their entropy, especially when arbitrary finitely generated groups are considered.

Including strengthening statement \ref{thmII}\ref{II:sft_ent}, there are likely other statements which could be added to Theorems \ref{thmI} and \ref{thmII}. The types of dynamical properties explored in this work are by no means exhaustive, so future work may be able to add to these theorems, and any such work will likely use free extensions extensively as they have been used here. In addition to extending these theorems, expanding the theory of free extensions may be fruitful in the study of shifts on groups. For instance, while the forward direction of Lemma \ref{FE:sft} is true in greater generality using the more general embeddings of Barbieri \cite{barbieri}, the reverse direction for the specific case of free extensions is a new result to the knowledge of the author. This, along with other results about free extensions also indicates that the study of SFTs on groups may be reduced to the study of SFTs on finitely generated groups for which every forbidden shape for the SFT is a generating set for the group. One might call such SFTs \textit{intrinsic} to the group, in the sense that they arise only when considering the entire group and not any of its subgroups.

Lastly, the mere existence of this theorem suggests that it may be possible to classify other dynamical properties by properties of the group. To the knowledge of the author, these results may be the only results in symbolic dynamics that gives implications about the group only from dynamical properties of the group, let alone a complete characterization of the group by dynamical properties. By assuming additional structure on the group, it may be possible to characterize other dynamical properties by this structure, and derive similar theorems as Theorems \ref{thmI} and \ref{thmII} for other classes of groups.

\subsection*{Acknowledgements}

The author would like to thank Robert Bland for fruitful conversations and Kevin McGoff for helpful recommendations about the structure and exposition of this paper. This research was done while the author was a PhD student at the University of North Carolina at Charlotte.

\bibliographystyle{abbrv}
\bibliography{refs}

\end{document}